# Operadic Definition of Non-Stricts Cells

Kamel Kachour

February 10, 2010

**Abstract**

In [11] we pursue Penon's work in higher dimensional categories by defining non-strict $\infty$-functors, non-strict natural $\infty$-transformations, and so on, all that with Penon's frameworks i.e with the "étirements catégoriques", where we have used an extension of this object, namely the "$n$-étirements catégoriques" ($n \in \mathbb{N}$). In this paper we are pursuing Batanin's work in higher dimensional categories by defining non-strict $\infty$-functors, non-strict natural $\infty$-transformations, and so on, using Batanin's frameworks i.e with the contractible operads, where we used an extension of this object, namely the globular colored contractible operads.

# Contents

The desire to understand non-strict cells also means the desire to go towards categorification within the theory of categories and, more widely, in mathematical categorification. Indeed, the main reason behind defining them comes from the following observation: the theory of categories is based on the notions of functors and natural transformations. The theory of the 2-categories is based on the notion of 2-functors, 2-natural transformations, and modifications, etc [5]. Therefore it is crucial to keep building on the higher categories by building on their higher functors, their higher natural transformations, as well as all of their higher cells. Thus the categorification begun by Baez and Dolan in [1] will be continued, as they had hoped. In [11], I prove that their construction is possible with an approach which is completely similar to the construction of Penon's non-strict $\infty$-categories. In this article we show that their construction is also possible with an approach that is completely similar to Batanin's construction of non-strict $\infty$-categories. Also this work probably can be done again for Leinster's non-strict $\infty$-categories. This work was exposed in Calais in June 2008 in the International Category Theory Conference, where in particular I introduced the important notion of contractible $T$-graphs ([10]) which are the extension of the classical ones [14]. The $T$-categories invented by Albert Burroni in 1971 organise all these constructions. Leinster then Hermida rediscovered these $T$-categories in 2004, [14] and 2000, [9]. Leinster calls them $T$-multicategories in [14]. The definition of Batanin's non-strict $\infty$-categories uses a globular non-symetrical alternative to May's operads. In this article the $T$-categories are built on the cartesian monad $T$ of strict $\infty$-category and are a globular, non-symetrical alternative to May's colored operads. Therefore they are the most suitable candidates so as to produce Batanin's non-strict cells and to carry on his work. The similar relation between operads and $T$-categories has been recently discovered in spite of the fact that both Albert Burroni and Peter May, Boardman and Vogt discovered $T$-categories and operads, respectively, at the same time [see 4, 8, 16]. Originally operads allow to give an algebraic structure to loop spaces and more recently many applications have been found, for example in Maxim Kontsevich's works [see 12, 13]. I would also like to mention Tom Leinster's work [14], which not only made Batanin's definition more accessible but also facilitated the rediscovery of $T$-categories. At the beginning of this paper his book will be referred to repeatedly because many of the defined concepts will be used here as well, which greatly condenses the text.

Michael Batanin built these non-strict $\infty$-categories with a contractible operad equipped with a composition system. I adopt the same point of view, using a countable infinity of globular colored operads, which will also be contractible and equipped with "bicolored composition systems" (called operation systems). Like Batanin, I have chosen to take them as initial objects, as this initiality happens to be crucial for constructing the sources and targets of the underlying graphs of the probable non-strict $\infty$-category of Batanin's non-strict $\infty$-categories. This construction proceeds in 4 stages: one first constructs a co-$\infty$-graph of operation systems, followed by a co-$\infty$-graph of globular colored operads, which will successively lead to a $\infty$-graph in the category of categories equipped with a monad, and finally the $\infty$-graph of algebras wished for, thanks to Eilenberg-Moore's functor. These algebras will contain all Batanin's non-strict cells: non-strict $\infty$-functors, non-strict $\infty$-natural transformations, non-strict $\infty$-modifications, etc.

The fundamental component of these constructions is the use of two colours: For example the operations system contains the symbols of operations of non-strict cells, and the action of monads of non-strict cells is then carried out on bicolored graphs. For instance, the operation system of non-strict 1-cells (i.e non-strict $\infty$-functors) contains twice the number of symbols of operations of non-strict $\infty$-categories, each one of them has a distinct colour. Furthermore, it contains the symbols of functors whose arity and co-arity have different colours, and finally it contains the bicolored operadic pointing. This bicoloring is obtained with the sum $T(1)+T(1)$ for globular arities and with the sum $1+1$ for co-arities, where 1 means the final $\infty$-graph. The promises of these constructions are kept for dimension 2. I demonstrate, in §6, that the weak 2-functors obtained are pseudo-2-functors and that the weak 2-natural-transformations are pseudo-2-natural transformations. Examples are left for the future.

One of the primary technical difficulties of this paper concerns the initial objects theorem (see theorem 1). This theorem is an adjunction result allowing the establishment of the existence of the free globular colored contractible operads of non-strict cells. We demonstrate this theorem in the same spirit as the proof suggested in Batanin [3, theorem 8.1 page 90], for the weakly initial contractible operad $K$, i.e "by doing the fusion" of two pairs of adjunction functors(see section 7.6). To do it we need to establish some results of com-pleteness and of cocompleteness of $T$-$\mathbb{C}\text{at}_c$ and $CT$-$\mathbb{G}\text{r}_{c,p}$ which are the two most important categories of this paper. The results are non-trivial, and if we want to resolve these difficulties we will need to develop the theory of Surcategories after the work of [7] (see section 7), where the surcategorical nature of these two categories enables us to draw conclusions about their completeness and cocompleteness. These "surcategorical" concepts are studied because the objects that we deal with in this paper are not fibrations (regrettably), and we need this notion to expand the work of Michael Batanin. Surcategories are the obvious generalisations of categories, and the "surcategorical theorems" that we establish in this paper (Freyd's Suradjoint theorem, Barr-Wells's Surcategorical theorem, and Beck's Surcategorical theorem), permit us to prove some properties for categories which have an underlying structure of surcategory (This can be done when these kinds of categories do not have enough properties to permit us to use the classic theorems on them). So theorem 1 is demonstrated in two very different stages: Firstly we construct two pairs of adjoints functors, which are in fact two pairs of suradjoints surfunctors (theorem 10 for the structure of monoid and theorem 11 for the structure of contractibility) and then establish a result (theorem 6), which is a generalization of techniques used in [3] which allows us to obtain from two pairs of "fusionnable" adjoints functors (and above all for two pairs of suradjoints surfunctors; see section 7.6), their "fusion", i.e to obtain in particular a single pair of adjoints functors which inherits from its respective structures (naturally the pairs of suradjoints surfunctors built in theorem 10 and in theorem 11 are fusionnable). The first pair of suradjoints surfunctors (section 8 and section 9) is more difficult to build than the second pair (section 10). We build it within the framework of the surmonoidals surcategories [7]. In [7] the authors establish an adjunction result (to obtain free "surmonoids") in an ideal context they label "numeral" [7, proposition 1.8 page 25]. We demonstrate a similar theorem (theorem 7) which also results from an ideal context that I label liberal and which allows us to establish a result of free surmonoids adapted to our particular situation. This theorem is an important piece in building the first pair of suradjoints surfunctors (theorem 10). The second pair of suradjoints surfunctors (theorem 11) is gullibly built by using elementary techniques deriving from logic. Let us indicate finally the attractive "fusion theorem" (theorem 5) which is a simple consequence of theorem 6, and which deserves to be studied more closely (although not used in this article) in the light of the very simple hypotheses if necessary (often realized in the practice) and of the "fusion result" which it proposes.

I am grateful to Jacques Penon who permitted me to access to the details of his conjoint work with Dominique Bourn ([7]). I am also grateful to Jean-Pierre Ledru and to my brother Ahmed Kachour for their support, and to Denis Bitouzé for his valuable advice in LaTeX. Finally I want to thanks my anonymous referee who not only help me to improve this paper but also who indicate me that it is well known that one can define $A_\infty$-functors between $A_\infty$-categories as algebras of appropriate contractible operads with two colours. I didn't know this result before and his or her remark convince me more that I found the good idea.

I dedicate this work to my uncle Mohamed-Kommandar Mezouar who taught us important life lessons even while struggling with a difficult illness.

# 1 Pointed and Contractibles $T$-Graphs

From here $\mathbb{T} = (T, \mu, \eta)$ refers to the cartesian monad of strict $\infty$-categories. Its cartesian feature permits us to build the bigategory $\text{Span}(T)$ of spans. The various concepts in this article are defined in this bicategory, which is described in Leinster [14, 4.2.1 page 138]. In all this paper if $\mathbb{C}$ is a category then $\mathbb{C}(0)$ is the class of its objects (but we often omit "0" when there is no confusion) and $\mathbb{C}(1)$ is the class of its morphisms. The symbol $:=$ means "by definition is".

## 1.1 $T$-Graphs

A $T$-graph $(C,d,c)$ is a datum of a diagram of $\infty$-$\mathbb{G}$r such as

$$T(G) \xleftarrow{d} C \xrightarrow{c} G$$

$T$-graphs are endomorphisms of Span$(T)$ and they form a category $T$-$\mathbb{G}$r (described in Leinster [14, definition 4.2.4 page 140]). If we choose $G \in \infty$-$\mathbb{G}$r$(0)$, the endomorphisms on $G$ (in Span$(T)$ ) forms a subcategory of $T$-$\mathbb{G}$r which will be noted $T$-$\mathbb{G}\mathrm{r}_G$, and it is well-known that $T$-$\mathbb{G}\mathrm{r}_G$ is a monoidal category such as the definition of its tensor:

$$(C,d,c)\bigotimes(C',d',c') := (T(C)\times_{T(G)} C', \mu(G)T(d)\pi_0, c\pi_1),$$

and its unity object $I(G) = (G, \eta(G), 1_G)$. We can remember that $I(G)$ is also an identity morphism of Span$(T)$ . We will see in section 9 other descriptions and useful properties of this tensorial product. The $\infty$-graph $G$ is called the graph of globular arities, and there is the arity functor $T$-$\mathbb{G}$r $\xrightarrow{\mathtt{A}}$ $\infty$-$\mathbb{G}$r, $(C,d,c) \longmapsto$ codomain$(c)$. We shall see furthermore that $(T$-$\mathbb{G}$r$, \mathtt{A})$ has several useful properties, as for example it is a liberal surmonoidal surcategory (see section 8). It is this property which allows us to obtain the first important monad of this article (see section 8.2 and section 9).

## 1.2 Pointed $T$-Graphs

A $T$-graph $(C,d,c)$ equipped with a morphism $I(G) \xrightarrow{p} (C,d,c)$ is called a pointed $T$-graph. Also we note $(C,d,c;p)$ for a pointed $T$-graph. That also means that one has a 2-cell $I(G) \xrightarrow{p} (C,d,c)$ of Span$(T)$ such as $dp = \eta(G)$ and $cp = 1_G$. We define in a natural way the category $T$-$\mathbb{G}\mathrm{r}_p$ of pointed $T$-graphs and the category $T$-$\mathbb{G}\mathrm{r}_{p,G}$ of $G$-pointed $T$-graphs: Their morphisms keep pointing in an obvious direction.

## 1.3 Contractible $T$-Graphs

Let $(C,d,c)$ be a $T$-graph. For any $k\in\mathbb{N}$ we consider

$$D_k=\{(\alpha,\beta)\in C(k)\times C(k)/s(\alpha)=s(\beta),t(\alpha)=t(\beta)\text{ and }d(\alpha)=d(\beta)\}$$

A contraction on that $T$-graph, is the datum, for all $k\in\mathbb{N}^*$, of a map

$$D_k\xrightarrow{[,]_k}C(k+1)$$

such that

- $s([\alpha,\beta]_k)=\alpha, t([\alpha,\beta]_k)=\beta$,
- $d([\alpha,\beta]_k)=1_{d(\alpha)=d(\beta)}$.

This maps $[,]_k$ form the bracket law (as the terminology in [11]). A $T$-graph which is equipped with a contraction will be called contractible and we note $(C,d,c;([,]_k)_{k\in\mathbb{N}^*})$ for a contractible $T$-graph. Nothing prevents a contractible $T$-graph from being equipped with several contractions. So here $CT\text{-}\mathbb{G}\mathrm{r}$ is a category of contractible $T$-graphs equipped with a specific contraction. The morphisms of this category preserves the contractions and one can also refer to the category $CT\text{-}\mathbb{G}\mathrm{r}_G$ where contractible $T$-graphs are only taken on a specific $\infty$-graph of globular arities $G$.

**Remark 1** If $(\alpha,\beta)\in D_k$ then this does not lead to $c(\alpha)=c(\beta)$, but this equality will be verified for constant $\infty$-graphs (see below) and in particular for collections with a finite number of colours in section 3.1 (These are the most important $T$-graphs in this article). We should also bear in mind $CT\text{-}\mathbb{G}\mathrm{r}_p$, the category of pointed and contractible $T$-graphs resulting from the previous definitions. A pointed and contractible $T$-graph will be noted $(C,d,c;([,]_k)_{k\in\mathbb{N}^*},p)$. □

### 1.4 Constant $\infty$-Graphs

A constant $\infty$-graph is a $\infty$-graph $G$ such as $\forall n,m \in \mathbb{N}$ we have $G(n) = G(m)$ and such as source and target maps are identity. We note $\infty$-$\mathbb{G}\mathrm{r}_c$ the corresponding category of constant $\infty$-graphs. Constant $\infty$-graph are important because it is in this context that we have an adjunction result (theorem 1) that we used to produce free colored contractibles operads of non-strict cells. We write $T$-$\mathbb{G}\mathrm{r}_c$ for the subcategory of $T$-$\mathbb{G}\mathrm{r}$ consisting of $T$-graphs with underlying $\infty$-graphs of globular arity which are constant $\infty$-graphs,$T$-$\mathbb{G}\mathrm{r}_{c,p}$ for the subcategory of $T$-$\mathbb{G}\mathrm{r}_p$ consisting of pointed $T$-graphs with underlying $\infty$-graphs of globular arity which are constant $\infty$-graphs, and we write $T$-$\mathbb{G}\mathrm{r}_{c,p,G}$ for the fiber subcategory in $T$-$\mathbb{G}\mathrm{r}_{c,p}$ (for a given $G$ in $\infty$-$\mathbb{G}\mathrm{r}_c$).

In proposition 17 we shall see that the pair $(T\text{-}\mathbb{G}\mathrm{r}_c, \mathtt{A})$ can be stuctured in a surmonoidal surcategory in [7] terminology.

# 2 Contractible $T$-Categories

### 2.1 $T$-Categories

A $T$-category is a monad of the bigategory $\mathrm{Span}(T)$ or in a equivalent way a monoid of the monoidal category $T$-$\mathbb{G}\mathrm{r}_G$ (for a specific $G$). The definition of $T$-categories are in Leinster [14, definition 4.2.2 page 140], and their category will be noted $T$-$\mathbb{C}\mathrm{at}$ and that of $T$-categories of the same $\infty$-graph of globular arities $G$ will be noted $T$-$\mathbb{C}\mathrm{at}_G$. A $T$-category $(B,d,c;\gamma,u) \in T$-$\mathbb{C}\mathrm{at}$ is specifically given by the morphism of (operadic) composition $(B,d,c)\bigotimes(B,d,c) \xrightarrow{\gamma} (B,d,c)$ and the (operadic) unit $I(G) \xrightarrow{u} (B,d,c)$ fitting axioms of associativity and unity [see 14]. Note that $(B,d,c;\gamma,u)$ has $(B,d,c;u)$ as natural underlying pointed $T$-graph.

## 2.2 Contractibles $T$-Categories and the Theorem of Initial Objects

A $T$-category $(B,d,c;\gamma,u)$ will be said to be contractible if its underlying $T$-graph is contractible. To specify the underlying contraction of contractible $T$-categories we eventually noted it $(B,d,c;\gamma,u,([,]_k)_{k\in\mathbb{N}^*})$. The category of contractible $T$-categories will be noted $CT$-$\mathbb{C}\mathrm{at}$, that of contractible $T$-categories of the same $\infty$-graph of globular arities $G$ will be noted $CT$-$\mathbb{C}\mathrm{at}_G$. We also write $CT$-$\mathbb{C}\mathrm{at}_c$ for the subcategory of $CT$-$\mathbb{C}\mathrm{at}$ whose objects are contractible $T$-categories whose underlying $\infty$-graph of globular arities is a constant $\infty$-graph. Besides there is an obvious forgetful functor

$$CT\text{-}\mathbb{C}\mathrm{at}_c \xrightarrow{O} T\text{-}\mathbb{G}\mathrm{r}_{c,p}$$

and there is the

**Theorem 1 (Theorem of Initial Objects)** *$O$ has a left adjoint $F$: $F \dashv O$. Moreover this left adjoint is in fact a left suradjoint.* □

See section 7 for the terminology adopted.

PROOF theorem 10 give the first surmonad $(L,\mathfrak{m},l)$, resulting from the suradjunction

$$(T\text{-}\mathbb{C}\mathrm{at}_c,\mathtt{A}) \underset{M}{\overset{U}{\rightleftarrows}} (T\text{-}\mathbb{G}\mathrm{r}_{c,p},\mathtt{A})$$

theorem 11 give the second surmonad $(C,m,c)$, resulting from the suradjunction

$$(CT\text{-}\mathbb{G}\mathrm{r}_{c,p},\mathtt{A}) \underset{H}{\overset{V}{\rightleftarrows}} (T\text{-}\mathbb{G}\mathrm{r}_{c,p},\mathtt{A})$$

The hypotheses of the theorem 6 are satisfied thanks to the following facts:

- Following theorem 10 and theorem 11 forgetful surfunctors $U$ and $V$ are surmonadic.

- It is easy to see that the surcategory $(T\text{-}\mathbb{G}\mathrm{r}_{c,p}, \mathtt{A})$ is surcomplet, surcocomplet, has $K$-equalizers, and that the surmonads $L$ and $C$ preserves $\vec{\mathbb{N}}$-surcolimits, thus by the theorem 3 the surcategories of Eilenberg-Moore suralgebras $(\mathbb{A}lg(L), \mathtt{A})$ and $(\mathbb{A}lg(C), \mathtt{A})$ have surcoequalizers and $\vec{\mathbb{N}}$-surcolimits. Thus surmonadicity of $U$ and $V$ implices that the surcategories $(T\text{-}\mathbb{C}\mathrm{at}_c, \mathtt{A})$ and $(CT\text{-}\mathbb{G}\mathrm{r}_{c,p}\mathtt{A})$ have surcoequalizers and $\vec{\mathbb{N}}$-surcolimits.

- It is easy to notice that the forgetful surfunctors $U$ and $V$ are faithfull and preserve $\vec{\mathbb{N}}$-surcolimits.

Thus this two suradjunctions are fusionable which permits, through theorem 6, to make the fusion

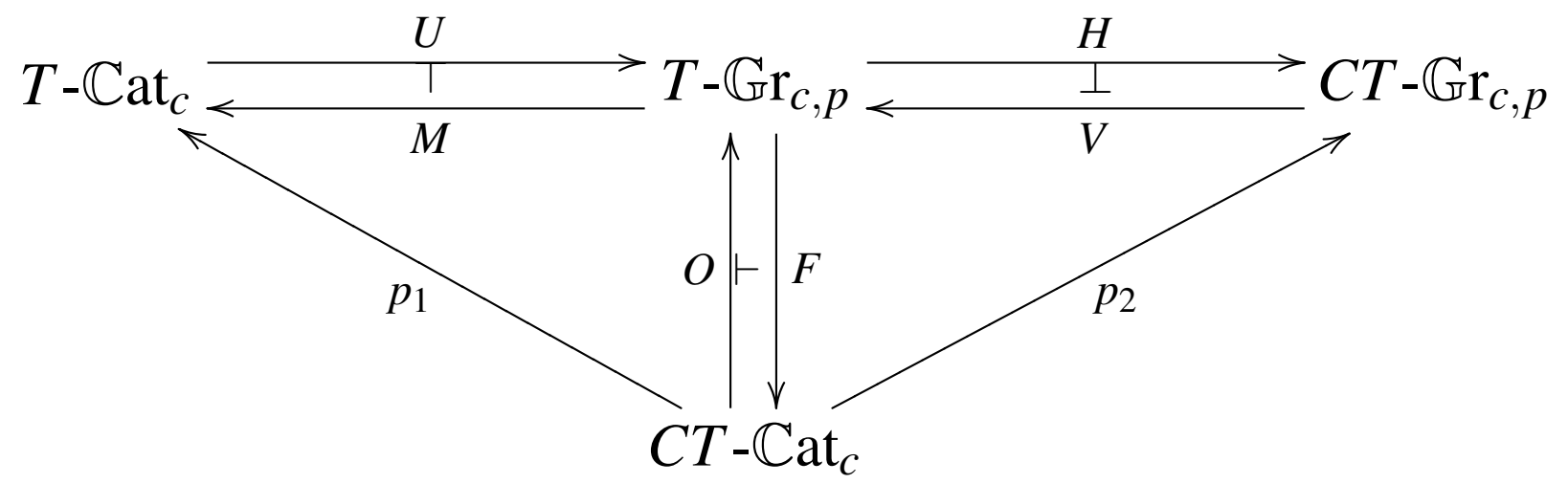

where trivially

$$CT\text{-}\mathbb{C}\mathrm{at}_c \simeq T\text{-}\mathbb{C}\mathrm{at}_c \times_{T\text{-}\mathbb{G}\mathrm{r}_{c,p}} CT\text{-}\mathbb{G}\mathrm{r}_{c,p}$$

■

The surmonad of this suradjunction $F \dashv O$ is noted $\mathbb{B} = (B, \rho, b)$.

**Remark 2** We can also prove that the forgetful functor

$$CT\text{-}\mathbb{C}\mathrm{at} \xrightarrow{O} T\text{-}\mathbb{G}\mathrm{r}_p$$

has a left adjoint but it is not a left suradjunction (see section 7). Indeed the problem is that this adjunction does not fix $\infty$-graphs of globular arity, and then does not fall within the framework of the Batanin's constructions where the pointed collection (composition system) and his free contractible operad has the same ∞-graphs of globular arity (who is 1, the final ∞-graph). As for the initial contractible operad of Batanin, we hope that our operation systems and their free colored contractible operads have the same ∞-graphs of globular arity.

To see that this left adjoint cannot fix ∞-graphs of globular arity we just need to take a $T$-graph $C$ whose ∞-graph of the globular arity $G$ is such as there is a $n \in \mathbb{N}$ where $G(p) = \emptyset$ if $p > n$. In that case a pair $(\alpha, \beta) \in D_n$ of $C$, will force the free $T$-category of $C$ to have an ∞-graph of globular arities $G'$ which is different from $G$, since $G'$ must contain $c([\alpha, \beta])$. □

### 2.3 $T$-Categories equipped with a System of Operations

Consider $(B, \overline{d}, \overline{c}; \gamma, u) \in T\text{-}\mathbb{C}\mathrm{at}_G$ and $(C, d, c) \in T\text{-}\mathbb{G}\mathrm{r}_G$. If there exists a diagram of $T\text{-}\mathbb{G}\mathrm{r}_G$

$$(I(G), \eta_G, id) \xrightarrow{p} (C, d, c) \xrightarrow{k} (B, \overline{d}, \overline{c})$$

such as $k \circ p = u$, then $(C, d, c)$ is qualified system of operations, and one can say that $(B, d, c; \gamma, u)$ is equipped with the system of operations $(C, d, c)$. With this definition and the previous theorem it is clear that all pointed $T$-graphs $(C, d, c; p)$ induces a free contractible $T$-category $F(C)$, which has $(C, d, c)$ as a system of operations. See also remark 7.

# 3 Version with $n$ Colours ($n \in \mathbb{N}$) of ∞-Graphs, Collections, Operads

### 3.1 $n$-Coloured ∞-Graphs ($n \in \mathbb{N}$)

For each $n \in \mathbb{N}$, we can possibly take the finite sum $1 \sqcup 1 ... 1 \sqcup 1$, $n$-times of the terminal ∞-graph 1. It will be said that this sum has $n$ colours, each colour being given by a ∞-graph 1 of the sum. The category of the $n$-coloured graphs (i.e. with $n$ colours) is the slice category $\infty$-$\mathbb{G}\mathrm{r}/1 \sqcup 1...1 \sqcup 1$. If $1 \overset{i_j}{\to} 1 \sqcup 1...1 \sqcup 1$ $(1 \leq j \leq n)$ refers to one of the canonical morphisms associated with this sum and if $G \overset{g}{\to} 1 \sqcup 1...1 \sqcup 1$ is an $n$-coloured $\infty$-graph, then the pullback $i_j^*(g)$ will refer to the underlying $\infty$-graph of $G$ with a $j$ colour. It will even be noted $i_j^*(G) = 1 \times_{1 \sqcup 1...1 \sqcup 1} G$. The functors $i_j^*$ are the colour functors.

## 3.2 Pointed $n$-Coloured Collections ($n \in \mathbb{N}$)

An $n$-coloured collection is a $T$-graph $(C,d,c)$ whose $\infty$-graph $G$ of globular arities is $1 \sqcup 1...1 \sqcup 1$. The category of the $n$-coloured collections will be noted $\mathbb{C}oll_{1,...,n} := T\text{-}\mathbb{G}\mathrm{r}_{1 \sqcup 1...1 \sqcup 1}$. If $(C,d,c) \in \mathbb{C}oll_{1,...,n}(0)$, the cells of $C$ are called symbols of globular operations. For $m \geq 0$ if $\alpha \in C(m)$, $d(\alpha)$ in a way is a "globular arity". Therefore the codomain of an operation is always monochromatic, i.e. always has a specific colour whereas its domain can be multicoloured. A pointed $n$-coloured collection is the datum of an object $(C,d,c) \in \mathbb{C}oll_{1,...,n}(0)$ and a morphism of $\mathbb{C}oll_{1,...,n}$

$$(1 \sqcup 1...1 \sqcup 1, \eta_{1 \sqcup 1...1 \sqcup 1}, id) \overset{p}{\to} (C,d,c).$$

The category of the pointed $n$-coloured collections is noted $\mathbb{C}oll_{1,...,n;p}$.

## 3.3 Contractible $n$-Coloured Operads ($n \in \mathbb{N}$)

An $n$-coloured operad is a monoid of $\mathbb{C}oll_{1,...,n}$. So an $n$-coloured operad $(B,d,c;\gamma,u)$ is composed of an underlying $n$-coloured collection $(B,d,c)$ and of two morphisms $u$, $\gamma$ of $\mathbb{C}oll_{1,...,n}$

$$(1 \sqcup ... \sqcup 1, \eta(1 \sqcup ... \sqcup 1), id) \overset{u}{\to} (B,d,c),$$
$$(T(B) \times_{T(1) \sqcup ... \sqcup T(1)} B, \mu(1 \sqcup ... \sqcup 1)T(d)\pi_0, c\pi_1) \overset{\gamma}{\to} (B,d,c)$$

fitting the axioms of monoids. Here $\mu$ refers to the multiplication of the monad $T$. If $(B',d',c';\gamma',u')$ is another $n$-coloured operad, a morphism

$$(B,d,c;\gamma,u) \xrightarrow{h} (B',d',c';\gamma',u')$$

is given by a morphism of $\infty$-graph $B \xrightarrow{h} B'$, such as we have the equalities $hu = u'$ and $h\gamma = \gamma'(h \bigotimes h)$. Here $h \bigotimes h$ refers to the single morphism induced by the universal pull-back property.

The category of the $n$-coloured operads is noted $\mathbb{O}per_{1\ldots,n}$ and that of contractible $n$-coloured operads is noted $C\mathbb{O}per_{1\ldots,n}$.

# 4 Systems of Operations of Non-Strict $n$-Cells

## 4.1 Preliminaries

The 2-coloured collection of the non-strict $n$-cells (for a given $n \in \mathbb{N}$) are just noted $C^n$ without specified its underlying structure, and we do the same simplification for its free contractible 2-coloured operads $B^n$.

From here on only the contractible 2-coloured operads of non-strict cells will be studied. All these operads are obtained applying the free functor of the theorem 1 to specific 2-coloured collections. These 2-coloured collections will be those of the non-strict cells and they count an infinite countable number of elements. Thus for each $n \in \mathbb{N}$ there is the 2-coloured collection of non-strict $n$-cells, $C^n$, which freely produces the free contractible 2-coloured operad $B^n$ of non-strict $n$-cells. The pointed collection $C^0$ is the system of composition of Batanin's operad of non-strict $\infty$-categories, i.e. the collection gathering all the symbols of atomic operations necessary for the non-strict $\infty$-categories, plus the symbols of operadic units (the latter are given by pointing). The pointed 2-coloured collection $C^1$ is adapted to non-strict $\infty$-functors, i.e. it gathers all the symbols of operations of the source and target non-strict $\infty$-categories (which will be composed of different colours whether they concern the source or the target). It also brings together the unary symbols of functors as well as the symbols of operadic units. Thus as we will see, the unary symbols of functors have a domain with the same colour as the domains and codomains of the symbols of operations of source non-strict $\infty$-categories and they have a codomain with the same colour as the domains and codomains of the symbols of operations of target non-strict $\infty$-categories. However these symbols of functors have domains and codomains with different colours. The pointed 2-coloured collection $C^2$ is adapted to natural non-strict $\infty$-transformations, etc.

## 4.2 Pointed 2-Coloured Collections $C^n(n \in \mathbb{N})$

In order to clearly see the bicolour feature of these symbols of operations, we write $(1+1)(n) := \{1(n), 2(n)\}$, which enables to identify $T(1) \sqcup T(1)$ with $T(1) \cup T(2)$ and $1 \sqcup 1$ with $1 \cup 2$. So the colour 1 and the colour 2 will be referred to. Let us move to the definition of $C^n(n \in \mathbb{N})$. In the diagram

$$T(1) \cup T(2) \xleftarrow{d} C^n \xrightarrow{c} 1 \cup 2$$

$C^n$ is a $\infty$-graph so that it contains both source and target applications which will be noted $C^n(m+1) \underset{t_m^{m+1}}{\overset{s_m^{m+1}}{\rightrightarrows}} C^n(m)$ , $(m \in \mathbb{N})$.

### 4.2.1 Definition of $C^0$

$C^0$ is Batanin's system of composition, i.e. there is the following collection $T(1) \overset{d^0}{\longleftarrow} C^0 \overset{c^0}{\longrightarrow} 1$ such as $C^0$ precisely contains the symbols of the compositions of non-strict $\infty$-categories $\mu_p^m \in C^0(m)(0 \leq p < m)$, plus the operadic unary symbols $u_m \in C^0(m)$. More specifically:

$\forall m \in \mathbb{N}$, $C^0$ contains the $m$-cell $u_m$ such as: $s_{m-1}^m(u_m) = t_{m-1}^m(u_m) = u_{m-1}$ (if $m \geq 1$); $d^0(u_m) = 1(m)(= \eta(1 \cup 2)(1(m)))$, $c^0(u_m) = 1(m)$.

$\forall m \in \mathbb{N} - \{0,1\}, \forall p \in \mathbb{N}$, such that $m > p$, $C^0$ contains the $m$-cell $\mu_p^m$ such as: If $p = m-1$, $s_{m-1}^m(\mu_{m-1}^m) = t_{m-1}^m(\mu_{m-1}^m) = u_{m-1}$. If $0 \leq p < m-1$, $s_{m-1}^m(\mu_p^m) = t_{m-1}^m(\mu_p^m) = \mu_p^{m-1}$. Also $d^0(\mu_p^m) = 1(m) \star_p^m 1(m)$, and inevitably $c^0(\mu_p^m) = 1(m)$.

Furthemore $C^0$ contains the 1-cell $\mu_0^1$ such as $s_0^1(\mu_0^1) = t_0^1(\mu_0^1) = u_0$, $d^0(\mu_0^1) = 1(1) \star_0^1 1(1)$, also inevitably $c^0(\mu_0^1) = 1(1)$.

The system of composition $C^0$ has got a well-known pointing $\lambda^0$ which is defined as $\forall m \in \mathbb{N}$, $\lambda^0(1(m)) = u_m$.

### 4.2.2 Definition of $C$

Firstly we will define a collection $(C,d,c)$ which will be useful to build the collections of non-strict cells. $C$ contains two copies of the symbols of $C^0$, each having a distinct colour: The symbols formed with the letters $\mu$ and $u$ are those of the colour 1, and those formed with the letters $\nu$ and $v$ are those of the colour 2. Let us be more precise:

$\forall m \in \mathbb{N}$, $C$ contains the $m$-cell $u_m$ such as: $s_{m-1}^m(u_m) = t_{m-1}^m(u_m) = u_{m-1}$ (if $m \geq 1$) and $d(u_m) = 1(m)$, $c(u_m) = 1(m)$.

$\forall m \in \mathbb{N} - \{0,1\}, \forall p \in \mathbb{N}$, such as $m > p$, $C$ contains the $m$-cell $\mu_p^m$ such as: If $p = m-1$, $s_{m-1}^m(\mu_{m-1}^m) = t_{m-1}^m(\mu_{m-1}^m) = u_{m-1}$. If $0 \leq p < m-1$, $s_{m-1}^m(\mu_p^m) = t_{m-1}^m(\mu_p^m) = \mu_p^{m-1}$. Also $d(\mu_p^m) = 1(m) \star_p^m 1(m)$, $c(\mu_p^m) = 1(m)$.

Furthemore $C$ contains the 1-cell $\mu_0^1$ such as $s_0^1(\mu_0^1) = t_0^1(\mu_0^1) = u_0$ and $d(\mu_0^1) = 1(1) \star_0^1 1(1)$, $c(\mu_0^1) = 1(1)$.

Besides, $\forall m \in \mathbb{N}$, $C$ contains the $m$-cellule $v_m$ such that: $s_{m-1}^m(v_m) = t_{m-1}^m(v_m) = v_{m-1}$ (if $m \geq 1$) and $d(v_m) = 2(m)$, $c(v_m) = 2(m)$.

$\forall m \in \mathbb{N} - \{0,1\}$, $\forall p \in \mathbb{N}$, such that $m > p$, $C$ contains the $m$-cell $v_p^m$ such as: If $p = m-1$, $s_{m-1}^m(v_{m-1}^m) = t_{m-1}^m(v_{m-1}^m) = v_{m-1}$. If $0 \leq p < m-1$, $s_{m-1}^m(v_p^m) = t_{m-1}^m(v_p^m) = v_p^{m-1}$. Also $d(v_p^m) = 2(m) \star_p^m 2(m)$, $c(v_p^m) = 2(m)$.

Furthemore $C$ contains the 1-cell $v_0^1$ such as $s_0^1(v_0^1) = t_0^1(v_0^1) = v_0$ and $d(v_0^1) = 2(1) \star_0^1 2(1)$, $c(v_0^1) = 2(1)$.

#### 4.2.3 Definition of $C^i (i = 1,2)$

$C^1$ is the system of operations of non-strict $\infty$-functors. It is built on the basis of $C$ adding to it a single symbol of functor(for each cell level):$\forall m \in \mathbb{N}$ the $F^m$ $m$-cell is added, which is such as: If $m \geq 1$, $s_{m-1}^m(F^m) = t_{m-1}^m(F^m) = F^{m-1}$. Also $d^1(F^m) = 1(m)$ and $c^1(F^m) = 2(m)$.

$C^2$ is the system of operations of non-strict natural $\infty$-transformations. $C^2$ is built on $C$, adding to it two symbols of functor (for each cell level) and a symbol of natural transformation. More precisely

$\forall m \in \mathbb{N}$ we add the $m$-cell $F^m$ such as: If $m \geq 1$, $s_{m-1}^m(F^m) = t_{m-1}^m(F^m) = F^{m-1}$. Also $d^2(F^m) = 1(m)$ and $c^2(F^m) = 2(m)$.

Then $\forall m \in \mathbb{N}$ we add the $m$-cell $H^m$ such as: If $m \geq 1$, $s_{m-1}^m(H^m) = t_{m-1}^m(H^m) = H^{m-1}$. Also $d^2(H^m) = 1(m)$ and $c^2(H^m) = 2(m)$.

And finally we add 1-cell $\tau$ such as: $s_0^1(\tau) = F^0$ and $t_0^1(\tau) = H^0$. Also $d^2(\tau) = 1_{1(0)}$ and $c^2(\tau) = 2(1)$.

We can point out that the 2-coloured collections $C^i$ $(i = 1,2)$ are naturally equipped with a pointing $\lambda^i$ defined by $\lambda^i(1(m)) = u_m$ and $\lambda^i(2(m)) = v_m$.

#### 4.2.4 Definition of $C^n$ for $n \geq 3$

In order to define the general theory of non-strict cells, it is necessary to define the systems of operations $C^n$ for the superior non-strict $n$-cells ($n \geq 3$).

This paragraph can be left out in the first reading. Each collection $C^n$ is built on $C$, adding to it the required cells. They contain four large groups of cells: The symbols of source and target non-strict $\infty$-categories, the symbols of operadic units (obtained on the basis of $C$), the symbols of functors (sources and targets), and the symbols of transformations (natural transformations, modification, etc). More precisely, on the basis of $C$:

**Symbols of Functors** $\forall m \in \mathbb{N}$, $C^n$ contains the $m$-cells $\alpha_0^m$ and $\beta_0^m$ such as: If $m \geq 1$, $s_{m-1}^m(\alpha_0^m) = t_{m-1}^m(\alpha_0^m) = \alpha_0^{m-1}$ and $s_{m-1}^m(\beta_0^m) = t_{m-1}^m(\beta_0^m) = \beta_0^{m-1}$. Furthermore $d^n(\alpha_0^m) = d^n(\beta_0^m) = 1(m)$ and $c^n(\alpha_0^m) = c^n(\beta_0^m) = 2(m)$.

**Symbols of Highers Cells (natural transformations, etc.)** $\forall p$, with $1 \leq p \leq n-1$, $C^n$ contains the $p$-cells $\alpha_p$ and $\beta_p$ which are such as: $\forall p$, with $2 \leq p \leq n-1$, $s_{p-1}^p(\alpha_p) = s_{p-1}^p(\beta_p) = \alpha_{p-1}$ and $t_{p-1}^p(\alpha_p) = t_{p-1}^p(\beta_p) = \beta_{p-1}$. If $p = 1$, $s_0^1(\alpha_1) = s_0^1(\beta_1) = \alpha_0^0$ and $t_0^1(\alpha_1) = t_0^1(\beta_1) = \beta_0^0$. What's more, $\forall p$, with $1 \leq p \leq n-1$, $d^n(\alpha_p) = d^n(\beta_p) = 1_p^0(1(0))$ and $c^n(\alpha_p) = c^n(\beta_p) = 2(p)$. Finally $C^n$ contains the $n$-cell $\xi_n$ such as $s_{n-1}^n(\xi_n) = \alpha_{n-1}$, $b_{n-1}^n(\xi_n) = \beta_{n-1}$ and $d^n(\xi_n) = 1_n^0(1(0))$ and $c^n(\xi_n) = 2(n)$ (Here $1_n^0$ is the map resulting from the reflexive structure of $T(1 \cup 2)$. See [11]).

We can see that $\forall n \in \mathbb{N}^*$, the 2-colored collection $C^n$ is naturally equipped with the pointing $1 \cup 2 \xrightarrow{\lambda^n} (C^n, d, c)$ defined as $\forall m \in \mathbb{N}, \lambda^n(1(m)) = u_m$ and $\lambda^n(2(m)) = v_m$.

## 4.3 The Co-$\infty$-Graph of Coloured Operads of Non-Strict Cells

In order not to make heavy notations we can write with the same notation $\delta_{n+1}^n$ and $\kappa_{n+1}^n$, sources and targets of the co-$\infty$-graph of coloured collections, the co-$\infty$-graph of coloured operads, and the $\infty$-graph in $\mathbb{M}nd$ below. There is no risk of confusion. The set $\{C^n/n \in \mathbb{N}\}$ has got a natural structure of co-$\infty$-graph. This co-$\infty$-graph is generated by diagrams

$$C^n \underset{\kappa^n_{n+1}}{\overset{\delta^n_{n+1}}{\rightrightarrows}} C^{n+1}$$

of pointed 2-coloured collections. For $n \geq 2$, these diagrams are defined as follows: First the $(n+1)$-colored collection contains the same symbols of operations as $C^n$ for $j$-cells, $0 \leq j \leq n-1$ or $n+2 \leq j < \infty$. For $n$-cells and $(n+1)$-cells the symbols of operations will change: $C^n$ contains the $n$-cell $\xi_n$ whereas $C^{n+1}$ contains $n$-cells $\alpha_n$ and $\beta_n$, in addition contains the $(n+1)$-cell $\xi_{n+1}$. If one notes $C^n - \xi_n$ the $n$-coloured collection obtained on the basis of $C^n$ by taking from it the $n$-cell $\xi_n$, then $\delta^n_{n+1}$ is defined as follows: $\delta^n_{n+1}|_{C^n-\xi_n}$ (i.e the restriction of $\delta^n_{n+1}$ to $C^n - \xi_n$) is the canonical injection $C^n - \xi_n \hookrightarrow C^{n+1}$ and $\delta^n_{n+1}(\xi_n) = \alpha_n$. In a similar way $\kappa^n_{n+1}$ is defined as follows: $\kappa^n_{n+1}|_{C^n-\xi_n} = \delta^n_{n+1}|_{C^n-\xi_n}$ and $\kappa^n_{n+1}(\xi_n) = \beta_n$. We can notice that $\delta^n_{n+1}$ and $\kappa^n_{n+1}$ keeps pointing, i.e we have for all $n \geq 1$ the equalities $\delta^n_{n+1}\lambda^n = \lambda^{n+1}$ and $\kappa^n_{n+1}\lambda^n = \lambda^{n+1}$.

The morphisms of 2-colored pointing collections of the diagram

$$C^O \underset{\kappa^0_1}{\overset{\delta^0_1}{\rightrightarrows}} C^1 \underset{\kappa^1_2}{\overset{\delta^1_2}{\rightrightarrows}} C^2 \underset{\kappa^2_3}{\overset{\delta^2_3}{\rightrightarrows}} C^3$$

have a similar definition:

By considering notation of section 4.2, we have for all integer $0 \leq p < n$ and for all $\forall m \in \mathbb{N}$:

$$\delta^0_1(\mu^n_p) = \mu^n_p;\ \delta^0_1(u_m) = u_m;\ \kappa^0_1(\mu^n_p) = v^n_p;\ \kappa^0_1(u_m) = v_m.$$

Also: $\delta^1_2(\mu^n_p) = \mu^n_p$; $\delta^1_2(u_m) = u_m$; $\delta^1_2(v^n_p) = v^n_p$; $\delta^1_2(v_m) = v_m$; $\delta^1_2(F^m) = F^m$. And $\kappa^1_2(\mu^n_p) = \mu^n_p$; $\kappa^1_2(u_m) = u_m$; $\kappa^1_2(v^n_p) = v^n_p$; $\kappa^1_2(v_m) = v_m$; $\kappa^1_2(F^m) = H^m$.

Finally: $\delta_3^2(\mu_p^n)=\mu_p^n$; $\delta_3^2(u_m)=u_m$; $\delta_3^2(\nu_p^n)=\nu_p^n$; $\delta_3^2(v_m)=v_m$; $\delta_3^2(F^m)=\alpha_0^m$; $\delta_3^2(H^m)=\beta_0^m$; $\delta_3^2(\tau)=\alpha_1$. And $\kappa_3^2(\mu_p^n)=\mu_p^n$; $\kappa_3^2(u_m)=u_m$; $\kappa_3^2(\nu_p^n)=\nu_p^n$; $\kappa_3^2(v_m)=v_m$; $\kappa_3^2(F^m)=\alpha_0^m$; $\kappa_3^2(H^m)=\beta_0^m$; $\kappa_3^2(\tau)=\beta_1$.

The pointed 2-coloured collections $C^n$ $(n\in\mathbb{N})$ are the sytems of operations of non-strict $n$-cells. Each of them freely produces the contractible 2-colored operads $B^n$ $(n\in\mathbb{N})$. Each of these contractible operads is equipped with a system of operations given by the pointed 2-coloured collection $C^n$. These operads $B^n$ are the operads of non-strict cells and are the most important objects in this article. They produce the monads $T_{B^n}$ whose algebras are the sought-after non-strict $n$-cells (see section 5 below). Due to the universal property of the unit $b$ of the monad $\mathbb{B}$, $C^n \xrightarrow{b(C^n)} B^n=B(C^n)$, one obtains the co-$\infty$-graph $B^\bullet$ of the coloured operads of non-strict cells.

$$
\begin{array}{ccccccccc}
B^0 & \overset{\delta_1^0}{\underset{\kappa_1^0}{\rightrightarrows}} & B^1 & \overset{\delta_2^1}{\underset{\kappa_2^1}{\rightrightarrows}} & B^2 & \cdots\rightrightarrows & B^{n-1} & \overset{\delta_n^{n-1}}{\underset{\kappa_n^{n-1}}{\rightrightarrows}} & B^n \cdots \\
{\scriptstyle b(C^0)}\uparrow & & \uparrow{\scriptstyle b(C^1)} & & \uparrow{\scriptstyle b(C^2)} & & \uparrow{\scriptstyle b(C^{n-1})} & & \uparrow{\scriptstyle b(C^n)} \\
C^0 & \overset{\delta_1^0}{\underset{\kappa_1^0}{\rightrightarrows}} & C^1 & \overset{\delta_2^1}{\underset{\kappa_2^1}{\rightrightarrows}} & C^2 & \cdots\rightrightarrows & C^{n-1} & \overset{\delta_n^{n-1}}{\underset{\kappa_n^{n-1}}{\rightrightarrows}} & C^n \cdots
\end{array}
$$

The commutativity property of these diagrams is important for the consistence of algebras (see section 5.5).

# 5 Monads and Algebras of Non-Strict Cells

$\mathbb{M}nd$ is the category of the categories equipped with a monad, and $\mathbb{A}dj$ is the category of the adjunction pairs. These categories are defined in [11].

## 5.1 Monads $T_{B^n}$ of Non Strict $n$-Cells $(n\in\mathbb{N})$.

Each $T$-category produces a monad which is described in [14, 4.3 page 150]. Hence $\forall n\in\mathbb{N}^*$, the operad $B^n$ of non-strict $n$-cells produce a monad $T_{B^n}$ on $\infty$-$\mathbb{G}$r$/1\cup 2$. More precisely, if we note $(B^n,d^n,c^n)$ its underlying $T$-graph we have (see notation used in section 9): $T_{B^n} := \Sigma_{c^n}(d^n)^*\widehat{T}$ (where we put $\widehat{T}(C,d,c) := (T(C),T(d),T(c))$). A bicolour $\infty$-graph $G \xrightarrow{g} 1\cup 2$ is often noted $G$ because there is no risk of confounding. We can therefore write $T_{B^n}(G)$ instead of $T_{B^n}(g)$, and it will be the same for the natural transformations $\delta_n^{n-1}$ and $\kappa_n^{n-1}$ (see below) and we write $T_{B^n}(G) := T(G)\times_{T(1\cup 2)} B^n$ (implied $T_{B^n}(g) = c^n\circ\pi_1$) and the definition of $T_{B^n}$ on morphisms is as easy. Projection on $T(G)\times_{T(1\cup 2)} B^n$ are noted $\pi_0$ and $\pi_1$. The definition of $T_{B^0}$ is similar.

## 5.2 The $\infty$-graph of $\mathbb{M}nd$ of Monads of Non-Strict Cells

Considering $G \xrightarrow{g} 1\cup 2$ , a bicolour $\infty$-graph. If we apply to it the monads $T_{B^n}$ and $T_{B^{n-1}}$ we obtain the equalities $d^n\pi_1 = T(g)\pi_0$, $d^{n-1}\pi_1 = T(g)\pi_0$. We also have $d^{n-1} = d^n\delta_n^{n-1}$ (To remove any confusion on our abuses of notations, the reader is encouraged to draw corresponding diagram). Thus we have $d^n\circ\delta_n^{n-1}\circ\pi_1 = 1_{T(1\cup 2)}\circ d^{n-1}\circ\pi_1 = 1_{T(1\cup 2)}\circ T(g)\circ\pi_0 = T(g)\circ 1_{T(G)}\circ\pi_0$. Hence the existence of a single morphism of $\infty$-graph

$$T(G)\times_{T(1\cup 2)} B^{n-1} \xrightarrow{\delta_n^{n-1}(G)} T(G)\times_{T(1\cup 2)} B^n$$

such as $\delta_n^{n-1}\pi_1 = \pi_1\delta_n^{n-1}(G)$ and $\pi_0 = \pi_0\delta_n^{n-1}(G)$. In particular we obtain the equality $c^n\pi_1\delta_n^{n-1}(G) = c^{n-1}\pi_1$. It is then easy to see that to each bicolour $\infty$-graph is associated the morphism(of $\infty$-$\mathbb{G}/1\cup 2$): $T_{B^{n-1}}(G) \xrightarrow{\delta_n^{n-1}(G)} T_{B^n}(G)$ (These morphisms are still simply called $\delta_n^{n-1}(G)$). It is very easy to see that the set of these morphisms produce a natural transformation $T_{B^{n-1}} \xrightarrow{\delta_n^{n-1}} T_{B^n}$. It is shown that $\delta_n^{n-1}$ fits the axioms $\mathbb{M}nd1$ and $\mathbb{M}nd2$ of the morphisms of monads (these axioms are in [11]) (particularly because $B^{n-1} \xrightarrow{\delta_n^{n-1}} B^n$ is a morphism of operads). Hence we can get the morphism of $\mathbb{M}nd$

$$(\infty\text{-}\mathbb{G}r/1\cup 2, T_{B^n}) \xrightarrow{\delta_n^{n-1}} (\infty\text{-}\mathbb{G}r/1\cup 2, T_{B^{n-1}})$$

Thus the morphisms of coloured operads $B^{n-1} \overset{\delta_n^{n1}}{\underset{\kappa_n^{n1}}{\rightrightarrows}} B^n$ $(n \geq 2)$, creates natural transformations $T_{B^{n-1}} \overset{\delta_n^{n-1}}{\underset{\kappa_n^{n1}}{\rightrightarrows}} T_{B^n}$ which fits into the axioms $\mathbb{M}nd1$ and $\mathbb{M}nd2$ of morphisms of monads. So we can get the diagrams of $\mathbb{M}nd(n \geq 2)$

$$(\infty\text{-}\mathbb{G}r/1 \cup 2, T_{B^n}) \overset{\delta_n^{n-1}}{\underset{\kappa_n^{n-1}}{\rightrightarrows}} (\infty\text{-}\mathbb{G}r/1 \cup 2, T_{B^{n-1}})$$

Similarly the morphisms $B^0 \overset{\delta_1^0}{\underset{\kappa_1^0}{\rightrightarrows}} B^1$ produce two natural transformations $T_{B^0} \circ i_1^* \overset{\delta_1^0}{\longrightarrow} i_1^* \circ T_{B^1}$, $T_{B^0} \circ i_2^* \overset{\kappa_1^0}{\longrightarrow} i_2^* \circ T_{B^1}$ ($i_1^*$ and $i_2^*$ are the colour functors) which also fits $\mathbb{M}nd1$ and $\mathbb{M}nd2$, which leads to the diagram of $\mathbb{M}nd$

$$(\infty\text{-}\mathbb{G}r/1 \cup 2, T_{B^1}) \overset{\delta_1^0}{\underset{\kappa_1^0}{\rightrightarrows}} (\infty\text{-}\mathbb{G}r/1 \cup 2, T_{B^0})$$

It is generally appeared that the building of the monad associated to a $T$-category is functorial, so the diagram of $\mathbb{M}nd$

$$\cdots \rightrightarrows (\infty\text{-}\mathbb{G}r/1 \cup 2, T_{B^n}) \cdots \rightrightarrows (\infty\text{-}\mathbb{G}r/1 \cup 2, T_{B^1}) \rightrightarrows (\infty\text{-}\mathbb{G}r/1 \cup 2, T_{B^0})$$

is a $\infty$-graph: The $\infty$-graph $T_{B^\bullet}$ of $\mathbb{M}nd$ of monads of non-strict cells.

## 5.3 The $\infty$-Graph of $\mathbb{C}AT$ of Batanin's Algebras of Non-Strict Cells

As in Kachour [11, § 4.3] we know that we have the functors

$$\mathbb{M}nd \xrightarrow{A} \mathbb{A}dj \xrightarrow{D} \mathbb{C}AT$$

where $A$ is the functor, which is linked with any monad, its pair of adjunction functors and where $D$ is the projection functor which associates $X$ with any adjunction $X \underset{F}{\overset{G}{\rightleftarrows}} Y$. So it is easy to see that $D \circ A$ associates its category of Eilenberg-Moore algebras to any monads. Particularly the functor $\mathbb{M}nd \xrightarrow{D\circ A} \mathbb{C}AT$ produces the following $\infty$-graph of $\mathbb{C}AT$

$$\cdots \mathbb{A}lg(T_{B^n}) \underset{\beta^n_{n-1}}{\overset{\sigma^n_{n-1}}{\rightrightarrows}} \mathbb{A}lg(T_{B^{n-1}}) \cdots \mathbb{A}lg(T_{B^1}) \underset{\beta^1_0}{\overset{\sigma^1_0}{\rightrightarrows}} \mathbb{A}lg(T_{B^0})$$

which is the $\infty$-graph $\mathbb{A}lg(T_{B^\bullet})$ of algebras of non-strict cells. It is the most important $\infty$-graph of this article since it contains all Batanin's non-strict cells.

## 5.4 Domains and Codomains of Algebras

Let us remember the morphisms of $\mathbb{M}nd$: $(C,T) \xrightarrow{(Q,t)} (C',T')$ are given by functors $C \xrightarrow{Q} C'$ and natural transformations $T' \circ Q \xrightarrow{t} Q \circ T$ whose fits $\mathbb{M}nd1$ and $\mathbb{M}nd2$. If we apply the functor $\mathbb{M}nd \xrightarrow{D\circ A} \mathbb{C}AT$ to these morphisms, one can get the functor, $\mathbb{A}lg(T) \to \mathbb{A}lg(T')$, defined on the objects as $(G,v) \longmapsto (Q(G), Q(v) \circ t(G))$. We can now describe the functors $\sigma^n_{n-1}$ and $\beta^n_{n-1}$ ($n \geq 1$):

- If $n \geq 2$ then $\mathbb{A}lg(T_{B^n}) \xrightarrow{\sigma^n_{n-1}} \mathbb{A}lg(T_{B^{n-1}})$, $(G,v) \longmapsto (G, v \circ \delta^{n-1}_n(G))$ and $\mathbb{A}lg(T_{B^n}) \xrightarrow{\beta^n_{n-1}} \mathbb{A}lg(T_{B^{n-1}})$, $(G,v) \longmapsto (G, v \circ \kappa^{n-1}_n(G))$.
- If $n = 1$ then $\mathbb{A}lg(T_{B^1}) \xrightarrow{\sigma^1_0} \mathbb{A}lg(T_{B^0})$, $(G,v) \longmapsto (i^*_1(G), i^*_1(v) \circ \delta^0_1(G))$ and $\mathbb{A}lg(T_{B^1}) \xrightarrow{\beta^1_0} \mathbb{A}lg(T_{B^0})$, $(G,v) \longmapsto (i^*_2(G), i^*_2(v) \circ \kappa^0_1(G))$.

## 5.5 Consistence of Algebras

As Penon's [1], Batanin's non-strict cells are particular in that they describe the hole semantics of their domain and codomain algebras as follows: If we have an algebra $(G,v)$ of non-strict $n$-cells, then a symbol of operation of the operad $B^n$ which has its counterpart in the operad $B^p$ ($0 \leq p < n$) will be semantically interpreted similarly via this algebra $(G,v)$ or via the algebra $\sigma_p^n(G,v)$ or the algebra $\beta_p^n(G,v)$.

More precisely, for all cells $(a,\alpha)$ of $T(G)\times_{T(1\cup 2)} B^p$ we have the following equalities:

$$v(\kappa_n^p(\alpha);a) = \beta_p^n(G,v)(\alpha;a) \text{ and } v(\delta_n^p(\alpha);a) = \sigma_p^n(G,v)(\alpha;a)$$

which is a simple consequence of the commutative property of the diagrams in section 4.3 applied to the bicolour $\infty$-graph $G$. To express this property, we will say that Batanin's algebras (as Penon's algebras) are consistent.

So as to illustrate this property, let us take for example the symbol of operation $H^m$ of the operad $B^2$ (identified with $b(C^2)(H^m)$). It will be semantically interpreted by an algebra $(G,v)\in \mathbb{A}lg(T_{B^2})$ on a $m$-cell $a\in G(m)$ (of colour 1), similarly to how the $F^m$ symbol of the $B^1$ operad is interpreted by the target algebra $\beta_1^2(G,v)\in \mathbb{A}lg(T_{B^1})$. Indeed the equalities $\kappa_2^1\pi_1 = \pi_1\kappa_2^1(G)$ and $\kappa_2^1 b(C^1) = b(C^2)\kappa_2^1$ immediately suggests that: $(a,F^m)\overset{\kappa_2^1(G)}{\longmapsto}(a,H^m)$, then $v(a,H^m) = (v\circ\kappa_2^1(G))(a,F^m) = \beta_1^2(G,v)(a,F^m)$, which expresses the property of consistence.

**Remark 3** This terminology is taken from measure theory where different coverings of a measurable subset are measured with the same value by a determined measure, which makes sense to that measure. □

# 6 Dimension 2

## 6.1 Dimension of Algebras

The dimension of Penon's algebras is defined in [17] and in [11]. The dimension of Batanin's algebras is totally similar, but we must precisely define the structures of the underlying $\infty$-magmas of these algebras so as to have a reflexive structure. So we can note $B^n \times_{T(1\cup 2)} T(G) \xrightarrow{v} G$ a $T_{B^n}$-algebra i.e a non-strict $n$-cell ($n \geq 1$). The two $\infty$-magmas ([11]) of this algebra are defined as follows: $\alpha \circ^n_p \beta := v(\mu^n_p; \eta(\alpha) \star^n_p \eta(\beta))$ and $1_\alpha := v([u_n, u_n]; 1_{\eta(\alpha)})$, if $\alpha, \beta \in G(n)$ and are with colour 1. Furthemore $\alpha \circ^n_p \beta := v(\nu^n_p; \eta(\alpha) \star^n_p \eta(\beta))$ and $1_\alpha := v([v_n, v_n]; 1_{\eta(\alpha)})$, if $\alpha, \beta \in G(n)$ and are with colour 2. Then $(G, v)$ has dimension 2 if its two underlying $\infty$-magmas has dimension 2. We have the same definition for $T_{B^0}$-algebras (i.e non-stricts $\infty$-categories).

## 6.2 The $T_{B^1}$-Algebras of dimension 2 are Pseudo-2-Functors

Let $(G, v)$ be a $T_{B^1}$-algebra of dimension 2. The $T_{B^0}$-algebra's source of $(G, v)$: $\sigma^1_0(G, v) = (i^*_1(G), i^*_1(v) \circ \delta^0_1(G))$ put on $i^*_1(G)$ a bicategory structure which coincides with the one produced by $(G, v)$ on $i^*_1(G)$. In the same way, the $T_{B^0}$-algebra target of $(G, v)$: $\beta^1_0(G, v) = (i^*_2(G), i^*_2(v) \circ \kappa^0_1(G))$ put on $i^*_2(G)$ a bicategory structure which coincides with that one produced by $(G, v)$ on $i^*_2(G)$. All these coincidences come from the consistence of algebras, and so we can therefore make all our calculations merely with the $T_{B^1}$-algebra $(G, v)$ to show the given below axiom of associativity-distributivity (that we call $AD$-axiom) of pseudo-2-functors. For other axioms of the pseudo-2-functors, which are easier, we proceed in the same way. Let $F^m (m \in \mathbb{N})$ be the unary operations symbols of functors of the operad $B^1$. The $T_{B^1}$-algebra of dimension 2 interprets these symbols into pseudo-2-functors. Indeed if $B^1 \times_{T(1\cup 2)} T(G) \xrightarrow{v} G$ is a $T_{B^1}$-algebra of dimension 2 then we get: $\forall m \in \mathbb{N}$, $F(a) := v(F^m; \eta(a))$ if $a \in G(m)$($a$ has the colour 1), which defines a morphism of $\infty$-graphs $i^*_1(G) \xrightarrow{F} i^*_2(G)$ where $i^*_1(G)$ and $i^*_2(G)$ are bicategories. So we will show that this morphism $F$ fits the $AD$-axiom of pseudo-2-functors. Let $x \xrightarrow{a} y \xrightarrow{b} z \xrightarrow{c} t$ be a 1-cellules diagram of $i^*_1(G)$, we are going to check that we get the following commutativity

$$\begin{array}{ccc}
F(a)\circ^1_0(F(b)\circ^1_0 F(c)) & \xRightarrow{1_{F(a)}\circ^2_0 d(b,c)} & F(a)\circ^1_0 F(b\circ^1_0 c)\\
{\scriptstyle a(F(a),F(b),F(c))}\Uparrow & & \Downarrow{\scriptstyle d(a,b\circ^1_0 c)}\\
(F(a)\circ^1_0 F(b))\circ^1_0 F(c) & \square & F(a\circ^1_0(b\circ^1_0 c))\\
{\scriptstyle d(a,b)\circ^2_0 1_{F(c)}}\Downarrow & & \Uparrow{\scriptstyle F(a(a,b,c))}\\
F(a\circ^1_0 b)\circ^1_0 F(c) & \xRightarrow{d(a\circ^1_0 b,c)} & F((a\circ^1_0 b)\circ^1_0 c)
\end{array}$$

where $a\circ^1_0(b\circ^1_0 c)\xRightarrow{a(a,b,c)}(a\circ^1_0 b)\circ^1_0 c$ is an associativity coherence cell and $F(a)\circ^1_0 F(b)\xRightarrow{d(a,b)} F(a\circ^1_0 b)$ is a distributivity coherence cell (particular to pseudo-2-functors). The strategy to demonstrate the $AD$-axiom is simple: We build a diagram of 3-cells of $B^1$ which will be semantically interpreted by the $T_{B^1}$-algebras of dimension 2 as the $AD$-axiom. To be clearer, the operadic multiplication of the coloured operad $B^1$

$$B^1\times_{T(1\cup 2)} T(B^1)\xrightarrow{\gamma} B^1$$

will be noted $\gamma_i$ for each $i$-cellular level. Let the following 2-cells in $B^1$:

$$d := [\gamma_1(\nu^1_0;\eta(F^1)\star^1_0\eta(F^1));\gamma_1(F^1;\eta(\mu^1_0))];$$

$$a_1 := [\gamma_1(\mu^1_0;\eta(\mu^1_0)\star^1_0\eta(u_1));\gamma_1(\mu^1_0;\,\eta(u_1)\star^1_0\eta(\mu^1_0))];$$

$$a_2 := [\gamma_1(\nu^1_0;\eta(\nu^1_0)\star^1_0\eta(v_1));\,\gamma_1(\nu^1_0;\eta(v_1)\star^1_0\eta(\nu^1_0))].$$

**Remark 4** The operation symbol $d$ is interpreted by the algebra as the distributivity coherence cells of the pseudo-2-functors. The symbols $a_1$ and $a_2$ are interpreted as the associativity coherence cells, the first one for the source $\infty$-category non-strict the second one for target $\infty$-category non-strict. □

Then we can consider the following 2-cells of $B^1$:

$$\rho_1 = \gamma_2(\nu_0^2;\eta([F^1;F^1]) \star_0^2 \eta(d));$$

$$\rho_2 = \gamma_2(d;1_{\eta(u_1)} \star_0^2 1_{\eta(\mu_0^1)});$$

$$\rho_3 = \gamma_2(F^2;\eta(a_1));$$

$$\rho_4 = \gamma_2(d;1_{\eta(\mu_0^1)} \star_0^2 1_{\eta(u_1)});$$

$$\rho_5 = \gamma_2(\nu_0^2;\eta(d) \star_0^2 \eta([F^1;F^1]))$$

$$\rho_6 = \gamma_2(a_2;1_{\eta(F^1)} \star_0^2 1_{\eta(F^1)} \star_0^2 1_{\eta(F^1)}).$$

This 2-cells are the conglomerations of operation symbols that are interpreted by algebras as the coherence 2-cells of the diagram of the *AD*-axiom of pseudo-2-functors

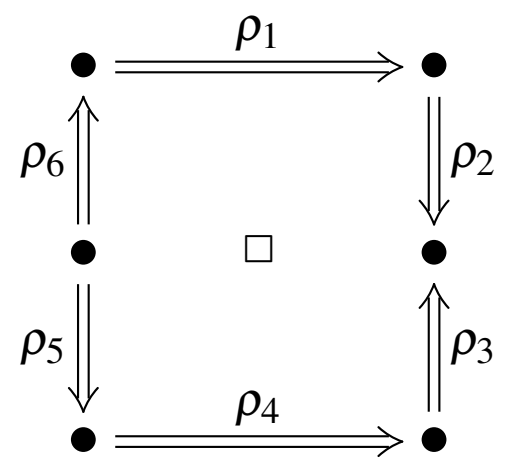

Then we consider the following 2-cells of $B^1$

$$\Lambda_1 = \gamma_2(\nu_1^2;\eta(\gamma_2(\nu_1^2;\eta(\rho_2) \star_1^2 \eta(\rho_1))) \star_1^2 \eta(\rho_6));$$

$$\Lambda_1' = \gamma_2(\nu_1^2;\eta(\rho_2) \star_1^2 \eta(\gamma_2(\nu_1^2;\eta(\rho_1) \star_1^2 \eta(\rho_6))));$$

$$\Lambda_2 = \gamma_2(\nu_1^2;\eta(\gamma_2(\nu_1^2;\eta(\rho_3) \star_1^2 \eta(\rho_4))) \star_1^2 \eta(\rho_5));$$

$$\Lambda_2' = \gamma_2(\nu_1^2;\eta(\rho_3) \star_1^2 \eta(\gamma_2(\nu_1^2;\eta(\rho_4) \star_1^2 \eta(\rho_5)))).$$

We can show that these 2-cells are parallels and with the same domain, so they are connected with coherences 3-cells

$\Theta_1 = [\Lambda_1, \Lambda_1']$, $\Theta_2 = [\Lambda_1', \Lambda_2]$, $\Theta_3 = [\Lambda_2, \Lambda_2']$,

and the interpretation by $T_{B^1}$-algebras of dimension 2 of this 3-cells gives the *AD*-axiom of pseudo-2-functors.

## 6.3 The $T_{B^2}$-Algebras of dimensions 2 are Natural Pseudo-2-Transformations

Let $(G,v)$ be a $T_{B^2}$-algebra of dimension 2. The $T_{B^0}$-algebra source of $(G,v)$: $\sigma_1^2(\sigma_0^1(G,v)) = (i_1^*(G), i_1^*(v\circ\delta_2^1(G))\circ\delta_1^0(G))$ put in $i_1^*(G)$ a bicategory structure which coincides with the one produced by $(G,v)$ on $i_1^*(G)$. In the same way, the $T_{B^0}$-algebra target of $(G,v)$: $\beta_1^2(\beta_0^1(G,v)) = (i_2^*(G), i_2^*(v\circ\kappa_2^1(G))\circ\kappa_1^0(G))$ put in $i_2^*(G)$ a bicategory structure which coincides with the one produced by $(G,v)$ on $i_2^*(G)$. And the $T_{B^1}$-algebra source of $(G,v)$: $\sigma_1^2(G,v) = (G, v\circ\delta_2^1(G))$ produces a pseudo-2-functor $F_1$ (it is the previous chapter) which coincides with the one produced by $(G,v)$ i.e the one built with the $\infty$-graph morphism $i_1^*(G)\xrightarrow{F_1} i_2^*(G)$ defined as: $F_1(a) := v(F^m;\eta(a))$ if $a\in i_1^*(G)(m)$. Besides the $T_{B^1}$-algebra target of $(G,v)$: $\beta_1^2(G,v) = (G, v\circ\kappa_2^1(G))$ produces a pseudo-2-functor $H_1$ which coincides with the one produced by $(G,v)$ i.e the one built with the $\infty$-graph morphism $i_1^*(G)\xrightarrow{H_1} i_2^*(G)$ defined as: $H_1(a) := v(H^m;\eta(a))$ if $a\in i_1^*(G)(m)$. All these coincidences come from the consistence of algebras, and we can therefore make all our calculations merely with the $T_{B^2}$-algebra $(G,v)$ (without using its source algebra or its target algebra) to show the axiom below of compatibility with associativity-distributivity of natural pseudo-2-transformations (that we call *CAD*-axiom). Then let $\tau$ be the unary operation symbol of natural transformation of the operad $B^2$. This symbol is interpreted by the $T_{B^2}$-algebras of dimension 2 as natural pseudo-2-transformations. Indeed if $B^2\times_{T(1\cup 2)} T(G)\xrightarrow{v} G$ is an $T_{B^2}$-algebra of dimension 2 then we write

$$\tau_1(a) := v(\tau; 1_{\eta(a)}),\ \text{if } a\in G(0)(a \text{ has colour} 1),$$

We can see that it defines a 1-cells family $\tau_1$ in $i_2^*(G)$ indexed by $i_1^*(G)(0)$

$$i_1^*(G) \underset{H_1}{\overset{F_1}{\rightrightarrows}} i_2^*(G) \quad (\Downarrow \tau_1)$$

We are going to show that the previous family $\tau_1$ fits the *CAD*-axiom of natural pseudo-2-transformations. For other axioms of natural pseudo-2-transformations, which are easier, we proceed in the same way. Let $x \xrightarrow{a} y \xrightarrow{b} z$ be an 1-cells diagram of $i_1^*(G)$, we are going to prove that we have the following commutativity

$$\begin{array}{ccc}
H_1(b)\circ_0^1(H_1(a)\circ_0^1\tau_1(x)) & \xRightarrow{1_{H_1(b)}\circ_0^2\omega(a)} & H_1(b)\circ_0^1(\tau_1(y)\circ_0^1 F_1(a)) \\
{\scriptstyle a(H_1(b),H_1(a),\tau_1(x))}\Big\Downarrow & & \Big\Downarrow{\scriptstyle a(H_1(b),\tau_1(y),F_1(a))} \\
(H_1(b)\circ_0^1 H_1(a))\circ_0^1\tau_1(x) & & (H_1(b)\circ_0^1\tau_1(y))\circ_0^1 F_1(a) \\
{\scriptstyle d_1(a,b)\circ_0^2 1_{\tau_1(x)}}\Big\Downarrow & \square & \Big\Downarrow{\scriptstyle \omega(b)\circ_0^2 1_{F_1(a)}} \\
H_1(b\circ_0^1 a)\circ_0^1\tau_1(x) & & (\tau_1(z)\circ_0^1 F_1(b))\circ_0^1 F_1(a) \\
{\scriptstyle \omega(b\circ_0^1 a)}\Big\Downarrow & & \Big\Downarrow{\scriptstyle a(\tau_1(z),F_1(b),F_1(a))} \\
\tau_1(z)\circ_0^1 F_1(b\circ_0^1 a) & \xLeftarrow[1_{\tau_1(z)}\circ_0^2 d_0(b,a)]{} & \tau_1(z)\circ_0^1(F_1(b)\circ_0^1 F_1(a)).
\end{array}$$

where in particular $H_1(a)\circ_0^1\tau_1(x) \xrightarrow{\omega(a)} \tau_1(y)\circ_0^1 F_1(a)$ is a coherence cell specific to natural pseudo-2-transformations. The strategy to demonstrate the *CAD*-axiom is similar to the previous demonstration (for the *AD*-axiom of pseudo-2-functors): We build a diagram of 3-cells of $B^2$ that will be semantically interpreted by the $T_{B^2}$-algebras of dimension 2 as the *CAD*-axiom. Like before operadic composition is

$$B^2 \times_{T(1\cup 2)} T(B^2) \xrightarrow{\gamma} B^2$$

will be noted $\gamma_i$ for each $i$-cellular level. So we can consider the following 2-cells of $B^2$

$$\omega := [\gamma_1(\nu_0^1;\eta(H^1)\star_0^1\eta(\tau));\gamma_1(\nu_0^1;\eta(\tau)\star_0^1\eta(F^1))];$$

$$d^F := [\gamma_1(\nu_0^1;\eta(F^1)\star_0^1\eta(F^1));\gamma_1(F^1;\eta(\mu_0^1))];$$

$$d^H := [\gamma_1(\nu_0^1;\eta(H^1)\star_0^1\eta(H^1));\gamma_1(H^1;\eta(\mu_0^1))];$$

$$a := [\gamma_1(\nu_0^1;\eta(\nu_1)\star_0^1\eta(\nu_0^1));\gamma_1(\nu_0^1;\eta(\nu_0^1)\star_0^1\eta(\nu_1))];$$

$$b := [\gamma_1(\nu_0^1;\eta(\nu_0^1)\star_0^1\eta(\nu_1));\gamma_1(\nu_0^1;\eta(\nu_1)\star_0^1\eta(\nu_0^1))].$$

Then we consider the following 2-cells

$$\rho_1 = \gamma_2(\nu_0^2;\eta([H^1;H^1])\star_0^2\eta(\omega));$$

$$\rho_2 = \gamma_2(a;1_{\eta(H^1)}\star_0^2 1_{\eta(\tau)}\star_0^2 1_{\eta(F^1)});$$

$$\rho_3 = \gamma_2(\nu_0^2;\eta(\omega)\star_0^2\eta([F^1;F^1]));$$

$$\rho_4 = \gamma_2(b;1_{\eta(\tau)}\star_0^2 1_{\eta(F^1)}\star_0^2 1_{\eta(F^1)});$$

$$\rho_5 = \gamma_2(\nu_0^2;\eta([\tau;\tau])\star_0^2\eta(d^F));$$

$$\rho_6 = \gamma_2(\omega;1_{\eta(\mu_0^1)});$$

$$\rho_7 = \gamma_2(\nu_0^2;\eta(d)\star_0^2\eta([\tau;\tau]));$$

$$\rho_8 = \gamma_2(a;1_{\eta(H^1)}\star_0^2 1_{\eta(H^1)}\star_0^2 1_{\eta(\tau)}).$$

We also consider one 2-cell $\rho_5'$ built as follows:

$$\delta^F := [\gamma_1(F^1;\eta(\mu_0^1));\gamma_1(\nu_0^1;\eta(F^1)\star_0^1\eta(F^1))].$$

In that case we define

$$\rho_5' = \gamma_2(\nu_0^2;\eta([\tau;\tau])\star_0^2\eta(\delta^F)).$$

These 2-cells are the conglomeration of operation symbols that are interpreted by algebras as the coherence 2-cells of the diagram of the *CAD*-axiom of natural pseudo-2-transformations

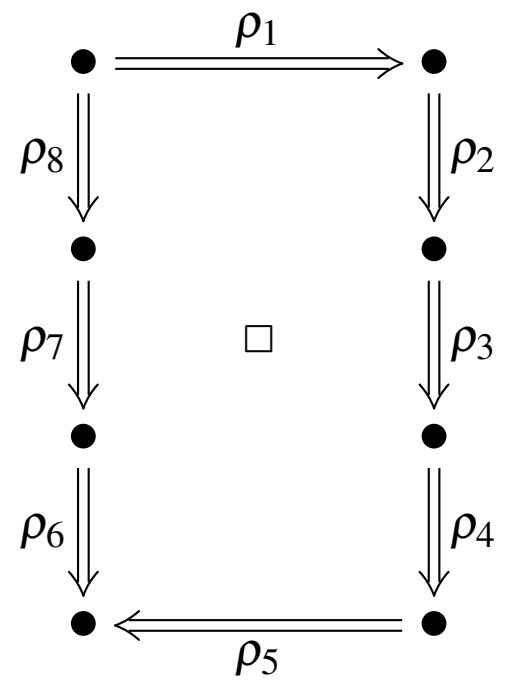

To built the ten coherence 2-cells $\Lambda_i (1 \leq i \leq 10)$ below, which enables to conclude, we need the following additional 2-cells

$$\Theta_1 = \gamma_2(\nu_1^2; \eta(\gamma_2(\nu_1^2; \eta(\nu_1^2) \star_1^2 \eta(\nu_2))) \star_1^2 \eta(\nu_2));$$

$$\Theta_2 = \gamma_2(\nu_1^2; \eta(\gamma_2(\mu_1^2; \eta(\nu_2) \star_1^2 \eta(\nu_1^2))) \star_1^2 \eta(\nu_2));$$

$$\Theta_3 = \gamma_2(\nu_1^2; \eta(\nu_2) \star_1^2 \eta(\gamma_2(\nu_1^2; \eta(\nu_1^2) \star_1^2 \eta(\nu_2))));$$

$$\Theta_4 = \gamma_2(\nu_1^2; \eta(\nu_2) \star_1^2 \eta(\gamma_2(\nu_1^2; \eta(\nu_2) \star_1^2 \eta(\nu_1^2))));$$

$$\Theta_5 = \gamma_2(\nu_1^2; \eta(\nu_1^2) \star_1^2 \eta(\nu_1^2)).$$

The 2-cells $\Lambda_i (1 \leq i \leq 10)$ are then defined in the following way

$$\Lambda_1 = \gamma_2(\Theta_1; \eta(\rho_4) \star_1^2 \eta(\rho_3) \star_1^2 \eta(\rho_2) \star_1^2 \eta(\rho_1));$$

$$\Lambda_2 = \gamma_2(\Theta_2; \eta(\rho_4) \star_1^2 \eta(\rho_3) \star_1^2 \eta(\rho_2) \star_1^2 \eta(\rho_1));$$

$$\Lambda_3 = \gamma_2(\Theta_3; \eta(\rho_4) \star_1^2 \eta(\rho_3) \star_1^2 \eta(\rho_2) \star_1^2 \eta(\rho_1));$$

$$\Lambda_4 = \gamma_2(\Theta_4; \eta(\rho_4) \star_1^2 \eta(\rho_3) \star_1^2 \eta(\rho_2) \star_1^2 \eta(\rho_1));$$

$$\Lambda_5 = \gamma_2(\Theta_5; \eta(\rho_4) \star_1^2 \eta(\rho_3) \star_1^2 \eta(\rho_2) \star_1^2 \eta(\rho_1)).$$

We can note as well $\lambda = \eta(\rho'_5) \star_1^2 \eta(\rho_6) \star_1^2 \eta(\rho_7) \star_1^2 \eta(\rho_8)$. And consider $\Lambda_6 = \gamma_2(\Theta_1;\lambda)$, $\Lambda_7 = \gamma_2(\Theta_2;\lambda)$; $\Lambda_8 = \gamma_2(\Theta_3;\lambda)$; $\Lambda_9 = \gamma_2(\Theta_4;\lambda)$; $\Lambda_{10} = \gamma_2(\Theta_5;\lambda)$.

We can prove that these 2-cells are parallels and with the same domain, so they are connected with coherences 3-cells: $\zeta_i := [\Lambda_i;\Lambda_{i+1}](1 \leq i \leq 9)$. And the interpretation by $T_{B^2}$-algebras of dimension 2 of these 3-cells gives the *CAD*-axiom of natural pseudo-2-transformations.

# 7 Theory of Surcategories

Theory of Surcategories is developped in [7] but many concepts are new, especially the notions of *K*-equalizers and *K*-coequalizers (see section 7.1 below) which allow to establish a surcategorical Freyd's Adjoint theorem (see section 7.3) which is a key result for the theorem 3 and thus for the theorem 1.

Let $\mathbb{G}$ a fixed category. A surcategory over $\mathbb{G}$ is an object of the 2-category $\mathbb{C}\mathrm{at}/\mathbb{G}$. Thus it is given by a couple $(\mathscr{C}, \mathtt{A})$, where $\mathscr{C}$ is a category and $\mathtt{A}$ is a functor (often called "arity functor", in reference to its use in this paper).In the whole continuation the arity functor is often noted with the letter $\mathtt{A}$ because there is no risk of confusion. The evident morphisms of $\mathbb{C}\mathrm{at}/\mathbb{G}$ are called surfunctors, but we also need in this paper morphim between surcategories with different base categories $\mathbb{G}$ and $\mathbb{G}'$. Therefore such a morphism $(\mathscr{C}, \mathtt{A}) \xrightarrow{(F,F_0)} (\mathscr{C}', \mathtt{A}')$ is given by two functors: $\mathscr{C} \xrightarrow{F} \mathscr{C}'$ and $\mathbb{G} \xrightarrow{F_0} \mathbb{G}'$ such as $\mathtt{A}'F = F_0\mathtt{A}$ (see for example section 7.3). For a fixed surcategory $(\mathscr{C}, \mathtt{A})$, its objects and its morphisms are respectively objects and morphisms of the underlying category $\mathscr{C}$.

The pairs of adjoints morphisms and monads in the 2-category $\mathbb{C}\mathrm{at}/\mathbb{G}$ are called respectively pairs of suradjoint surfunctors and surmonads. In fact every surcategorical concept will be often writing by adjoining "sur" before the categorical concept that it expands. But we sometimes forget the word "sur", when the context implies that no confusion is possible. It is easy to see that the category of algebras for a given surmonad is a surcategory. The objects of this surcategory will be named suralgebras.

We are going to see that most of the notions in the 2-category $\mathbb{C}$at will be done again in the 2-category $\mathbb{C}\mathrm{at}/\mathbb{G}$, and it is very likely that most of concepts and theorems in $\mathbb{C}$at extend in $\mathbb{C}\mathrm{at}/\mathbb{G}$. We will particularly demonstrate three theorems in $\mathbb{C}\mathrm{at}/\mathbb{G}$ coming from three important theorems in $\mathbb{C}$at: Freyd's Suradjoint theorem (which is the surcategorical version of Freyd's Adjoint classical theorem. See theorem 2), Barr-Wells's Surcategorical Theorem (which is the surcategorical version of the result that we can find in Barr and Wells [2, § 9.3]. See theorem 3), and Beck's Surcategorical Theorem (which is the surcategorical version of Beck's classical theorem. See theorem 4). These theorems are the generalisations of the classical ones.

## 7.1 Definition of Sur(co)limits

In Bourn and Penon [7, 1.3.2 page 23] the notions of limits and colimits in $\mathbb{C}\mathrm{at}/\mathbb{G}$ are defined and these notions will be used afterwards. To facilitate the reader we will remind the definitions.

If $\mathbb{C}$ is a small category and if $\mathscr{E}$ is a category, then we have the classical diagonal functor $\mathscr{E} \xrightarrow{\Delta} \mathscr{E}^{(\mathbb{C})}$, which sends an object to a constant functor and which sends a morphism to a constant natural transformation.

Moreover if $(\mathscr{E}, \mathtt{A})$ is a surcategory, let $\mathscr{E}^{(\mathbb{C})}$ be the subcategory of $\mathscr{E}^{\mathbb{C}} \times \mathbb{G}$ given by:

$$\mathscr{E}^{(\mathbb{C})}(0) = \{(F,B) \in \mathscr{E}^{\mathbb{C}} \times \mathbb{G} / \mathtt{A}F = \Delta(B)\},$$

$$\mathscr{E}^{(\mathbb{C})}(1) = \{(F,B) \xrightarrow{(\tau,b)} (F',B')/b \in \mathbb{G}(1) \text{ and } \tau \text{ is a natural transformation such as } \mathtt{A}\tau = \Delta(b)\}.$$

$\mathscr{E}^{(\mathbb{C})}$ has a natural surcategory structure given by the second projection: $\mathscr{E}^{(\mathbb{C})} \xrightarrow{\mathtt{A}} \mathbb{G}$, $(F,B) \longmapsto B$. In fact $(\mathscr{E}^{(\mathbb{C})}, \mathtt{A})$ is a cotensor of the $\mathbb{C}$at-enriched category $\mathbb{C}\mathrm{at}/\mathbb{G}$ ($\mathbb{C}\mathrm{at}/\mathbb{G}$ is a $\mathbb{C}\mathrm{at}$-enriched category because it is a 2-category), because we have the following isomorphism in $\mathbb{C}\mathrm{at}$

$$\mathbb{C}\mathrm{at}/\mathbb{G}((\mathscr{E}',\mathtt{A});(\mathscr{E}^{(\mathbb{C})},\mathtt{A})) \simeq Funct(\mathbb{C};\mathbb{C}\mathrm{at}/\mathbb{G}((\mathscr{E}',\mathtt{A});(\mathscr{E},\mathtt{A}))).$$

We also have the diagonal surfunctor (also noted $\Delta$): $(\mathscr{E},\mathtt{A}) \xrightarrow{\Delta} (\mathscr{E}^{(\mathbb{C})},\mathtt{A})$ defined by $x \longmapsto (\Delta(x),\mathtt{A}(x))$. If $(F,B) \in (\mathscr{E}^{(\mathbb{C})},\mathtt{A})$, a surcone of $(F,B)$ is a morphism $\Delta(x) \xrightarrow{(\tau,b)} (F,B)$ ($x \in \mathscr{E}$) of $(\mathscr{E}^{(\mathbb{C})},\mathtt{A})$, where $\mathtt{A}(x) \xrightarrow{b} B$ is a morphism of $\mathbb{G}$. In the same way we define surcocones.

It is easy to see that if $\mathbb{C}$ is connected, then every surcone is a cone in the classical way (respectively, every surcocone is a cocone in the classical way).

The surcategory $(\mathscr{E},\mathtt{A})$ has got $\mathbb{C}$-surlimits (that [7] call $(\mathbb{C})$-limits) if every $(F,B) \in (\mathscr{E}^{(\mathbb{C})},\mathtt{A})$ has a universal surcone $\Delta(x) \xrightarrow{(\tau,1_B)} (F,B)$ such as $\mathtt{A}(x) = B$, i.e if we give ourselves an other surcone $\Delta(y) \xrightarrow{(\sigma,b)} (F,B)$ (where $\mathtt{A}(y) \xrightarrow{b} B$ is a morphism of $\mathbb{G}$) then there is a unique morphism $y \xrightarrow{f} x$ in $\mathscr{E}$ such that $(\tau,1_B)\Delta(f) = (\sigma,b)$. The definition of $\mathbb{C}$-surcolimits are duals. The definition of surlimits and surcolimits enable us to include the case where $\mathbb{C}$ is the empty category, which allows to give an alternative definition of the surinitial objects (see section 7.3) and the surfinal objects.

If $\mathbb{C}$ is connected and nonempty then it is easy to note that the following definitions are equivalent

- $(\mathscr{E},\mathtt{A})$ has $\mathbb{C}$-surlimits.
- $\forall(F,B) \in (\mathscr{E}^{(\mathbb{C})},\mathtt{A})$, $(F,B)$ has an universal surcone $\Delta(x) \xrightarrow{(\tau,1_B)} (F,B)$ such as $\mathtt{A}(x) = B$.
- $\forall(F,B) \in (\mathscr{E}^{(\mathbb{C})},\mathtt{A})$,the functor $\mathbb{C} \xrightarrow{F} \mathscr{E}_B$ has a limit which is preserved by the canonical inclusion $\mathscr{E}_B \hookrightarrow \mathscr{E}$.
- The diagonal surfunctor $(\mathscr{E},\mathtt{A}) \xrightarrow{\Delta} (\mathscr{E}^{(\mathbb{C})},\mathtt{A})$ has a right suradjoint.

In the same way, if $\mathbb{C}$ is connected and nonempty we have dual definitions for $\mathbb{C}$-surcolimits.

**Remark 5** Let $\vec{\mathbb{N}}$ be the category of the non-negative integers with the natural order. In the terminology adopted in [7] $\vec{\mathbb{N}}$-limits are colimits. We prefer to adopt the word $\vec{\mathbb{N}}$-colimits for this specific filtered colimit. And in the surcategorical context we prefer the word $\vec{\mathbb{N}}$-surcolimits instead of $(\vec{\mathbb{N}})$-colimits (as it is adopted by [7]). □

We are now going to define $K$-equalizers and $K$-coequalizers which are important notions because with them we get a suradjonction result similar to Freyd's Adjoint theorem (theorem 2), but more general. In particular the surcategory $(T\text{-}\mathbb{G}\mathrm{r}, \mathtt{A})$ trivially has $K$-equalizers and $K$-coequalizers (like every subsurcategory of it which are important for our study. See section 1.1, section 1.2, and section 1.3).

A surcategory $(\mathscr{E}, \mathtt{A})$ has $K$-equalizers if every pair $a \overset{f}{\underset{g}{\rightrightarrows}} b$ , which has the property $\mathtt{A}(f) = \mathtt{A}(g)$, has a equalizer $e$ in $\mathscr{E}$

$$\begin{array}{c} c \\ \downarrow e \\ a \overset{f}{\underset{g}{\rightrightarrows}} b \end{array}$$

such as $\mathtt{A}(e) = \mathtt{A}(1_a)$. The definition of $K$-coequalizers is dual.

If $T$ is a surmonad on $(\mathscr{E}, \mathtt{A})$, Eilenberg-Moore algebras category $\mathscr{E}^T$ is trivially a surcategory $(\mathscr{E}^T, \mathtt{A})$ where its objects are called suralgebras, not only to precise the surcategorical context, but also to focus on the fact that a suralgebra is an algebra which lives in a fiber.

The following propositions are immediate and do not thus require detailed proof.

**Proposition 1** *Lets call split surfork, a split fork in the surcategorical context, i.e it is a diagram* $a \overset{f}{\underset{g}{\rightrightarrows}} b \xrightarrow{h} c$ *which is a fork in a fiber* $\mathscr{E}_B$. *Then such split surforks are absolute surcoequalizers.* □

**Proposition 2** *Every suralgebra (for a fixed surmonad) is a surcoequalizer.* □

**Proposition 3** $(\mathscr{E}, \mathtt{A})$ *is surcomplet iff* $(\mathscr{E}, \mathtt{A})$ *has surequalizers and surproducts.* □

**Proposition 4** $(\mathscr{E}, \mathtt{A})$ *is surcomplet iff* $(\mathscr{E}, \mathtt{A})$ *has surcoequalizers and sursums.* □

## 7.2 Some results of Sur(co)completeness of Suralgebras

The following propositions are very similar to the classical ones and do not thus require detailed proof.

**Proposition 5** *Let $T$ be a surmonad on $(\mathscr{E}, \mathtt{A})$. In this case:*

$$(\mathscr{E}, \mathtt{A}) \text{ is surcomplet} \Longrightarrow (\mathscr{E}^T, \mathtt{A}) \text{ is surcomplet}$$ □

**Proposition 6** *Let $T$ be a surmonad on $(\mathscr{E}, \mathtt{A})$. In this case:*

$$(\mathscr{E}, \mathtt{A}) \text{ has } K\text{-equalizers} \Longrightarrow (\mathscr{E}^T, \mathtt{A}) \text{ has } K\text{-equalizers}$$ □

**Proposition 7** *Let $T$ be a surmonad on $(\mathscr{E}, \mathtt{A})$. We suppose that $(\mathscr{E}, \mathtt{A})$ is surcocomplet. In this case:*

$$(\mathscr{E}^T, \mathtt{A}) \text{ has surcoequalizers} \Longleftrightarrow (\mathscr{E}^T, \mathtt{A}) \text{ is surcocomplet}$$ □

## 7.3 Freyd's Adjoint Theorem in the Surcategorical Context

As we are going to see, Freyd's Adjoint theorem remains true in the context of $\mathbb{C}\mathrm{at}/\mathbb{G}$. We call it "Freyd's Suradjoint theorem" to refer to its surcategorical nature. It is used for the proof of the theorem 3 which permits us to prove some surcocompleteness results that are important to apply the "Surcategorical Fusion Theorem" (see section 7.6) for the two important suradjunctions of this paper (see theorem 10 and theorem 11). But as we will demonstrate, unlike "Beck's Theorem in the Surcategorical Context" (see section 7.5), Freyd's Suradjoint theorem requires in addition $K$-equalizers (see theorem 2).

Let $(\mathscr{A},\mathtt{A}) \xrightarrow{F} (\mathscr{B},\mathtt{A})$ a surfunctor and $B \in (\mathscr{B},\mathtt{A})$. An object of the comma category $(B \downarrow F)$ is given by a couple $(A,a)$ which corresponds to a morphism $B \xrightarrow{a} F(A)$ in $\mathscr{B}$ and a morphism of $(B \downarrow F)$ is given by an arrow $(A,a) \xrightarrow{f} (A',a')$ such as $F(f)a = a'$.

The comma category $(B \downarrow F)$ is a surcategory. Indeed we have the arity functor $(B \downarrow F) \xrightarrow{\mathtt{A}} \mathtt{A}(B)/\mathbb{G}$ defined on the objects as: $(A,a) \longmapsto \mathtt{A}(a)$ and defined on the morphism as: $f \longmapsto \mathtt{A}(f)$ ($\mathtt{A}$ is here the arity functor of the surcategory $(\mathscr{B},\mathtt{A})$).

Furthermore we have the following canonical morphism of surcategories, given by the first projection

$$\begin{array}{ccc} (B \downarrow F) & \xrightarrow{Q} & \mathscr{A} \\ \downarrow{\scriptstyle \mathtt{A}} & & \downarrow{\scriptstyle \mathtt{A}} \\ \mathtt{A}(B)/\mathbb{G} & \xrightarrow{Q_0} & \mathbb{G} \end{array}$$

**Proposition 8** *Let* $(\mathscr{A},\mathtt{A}) \xrightarrow{G} (\mathscr{X},\mathtt{A})$ *be a surfunctor such as* $(\mathscr{A},\mathtt{A})$ *is surcomplet and has K-equalizers. We suppose that G preserves surlimits and K-equalizers. Then* $\forall B \in \mathscr{X}$*, the comma surcategory* $((B \downarrow G),\mathtt{A})$ *is surcomplet and has K-equalizers.* □

PROOF It is enough to prove that the functor $((B \downarrow G),\mathtt{A}) \xrightarrow{Q} (\mathscr{A},\mathtt{A})$ creates small surproducts, surequalizers, and $K$-equalizers. First we consider all functors $J \xrightarrow{F} (B \downarrow G)$ such as $F \in (B \downarrow G)^{(J)}$. Thus $QF \in A^{(J)}$ and if $J$ is a small discret category, then $limQF$ exists because $(\mathscr{A},\mathtt{A})$ is surcomplet. It is easy to prove (as in [15]) that $limF$ exists and that it is unique such as $Q(\lim F) = \lim QF$. If $J = \downdownarrows$ so we use a similar argument to prove that $Q$ creates surequalizers.

To prove that $Q$ creates $K$-equalizers we use a similar argument, but we must take $J = \downdownarrows$ and $F$ such as the image of the functor $\mathtt{A}F$ is a fixed arrow in $\mathtt{A}(B)/\mathbb{G}$. ■

Let $(\mathscr{D}, \mathtt{A})$ be a surcategory and let $G \in \mathbb{G}$. The object $0_G \in \mathscr{D}_G$ is surinitial if for all objects $d \in \mathscr{D}$, and for all $G \xrightarrow{b} \mathtt{A}(d)$ in $\mathbb{G}(1)$, there is a unique morphism $0_G \xrightarrow{x} d$ of $\mathscr{D}$ over $b$.

**Proposition 9** *Let* $(\mathscr{A}, \mathtt{A}) \xrightarrow{F} (\mathscr{B}, \mathtt{A})$ *be a surfunctor,* $B \in (\mathscr{B}, \mathtt{A})$*, and* $(R_B, v)$ *is an object of* $((B \downarrow F), \mathtt{A})$ *such as* $\mathtt{A}(v) = 1_{\mathtt{A}(B)}$*. In this case:*

$$(R_B, v) \text{ is surinitial in } ((B \downarrow F), \mathtt{A}) \iff v \text{ is initial in } (B \downarrow F) \qquad \square$$

**Lemma 1 (Lemma of the surinitial object)** *Let* $(\mathscr{D}, \mathtt{A})$ *a surcategory surcomplet with K-equalizers, and let* $G \in \mathbb{G}$*.*

*In this case we have the following equivalence*

| | | | |
|---|---|---|---|
| $(\mathscr{D}, \mathtt{A})$ *have a surinitial object in one fiber* $\mathscr{D}_G$ | $\iff$ | *There is a set I and a family of objects* $k_i \in \mathscr{D}_G$ *(*$i \in I$*) such as* $\forall d$ *in* $(\mathscr{D}, \mathtt{A})$*,* $\forall G \xrightarrow{h} \mathtt{A}(d)$ *in* $\mathbb{G}$*, there is an* $i \in I$*, there is a morphism* $k_i \to d$ *in* $\mathscr{D}$ *over h (via the arity functor).* | $\square$ |

The proof of this lemma is very similar to the classical one (see [7, proposition 1.8 page 25]) and thus it is not necessary to give the details of the demonstration. It is useful to note that this demonstration requires $K$-equalizers.

Let $(\mathscr{A}, \mathtt{A}) \xrightarrow{F} (\mathscr{B}, \mathtt{A})$ a surfunctor. An object $B \in (\mathscr{B}, \mathtt{A})$ have a solution set condition for $F$ if there is a set $I$ and a set of objects $\{(A_i, b_i)/i \in I \text{ and } \mathtt{A}(b_i) = 1_{\mathtt{A}(B)}\} \subset (B \downarrow F)$, such that $\forall (A, b) \in (B \downarrow F)$, $\exists i \in I$, $\exists A_i \xrightarrow{a_i} A$ in $(\mathscr{A}, \mathtt{A})$, such as $F(a_i)b_i = b$.

**Theorem 2 (Freyd's Suradjoint theorem)** *Let* $(\mathscr{A},\mathtt{A})$ *a surcomplet surcategory with* $K$*-equalizers, and let* $(\mathscr{A},\mathtt{A}) \xrightarrow{F} (\mathscr{B},\mathtt{A})$ *a surfunctor. In that case the following properties are equivalent*

| | | | |
|---|---|---|---|
| *F have a left suradjunction* | $\Longleftrightarrow$ | *F preserve surlimits and K-equalizers and every object* $B \in (\mathscr{B},\mathtt{A})$ *have a solution set condition for F* | □ |

PROOF First we suppose that $F$ preserves surlimits and $K$-equalizers and every object $B \in (\mathscr{B},\mathtt{A})$ has a solution set condition for $F$. Let $B \in Ob(\mathscr{B})$, the surcategory $(\mathscr{A},\mathtt{A})$ is surcomplet and have $K$-equalizers which are preserved by $F$, thus thanks to the proposition 8 we know that $((B \downarrow F),\mathtt{A})$ is surcomplet and have $K$-equalizers. Therefore $((B \downarrow F),\mathtt{A})$ verifies in addition the hypothesis "solution set condition" of the lemma of the surinitial object in the fiber $(B \downarrow F)_{1_{\mathtt{A}(B)}}$. Thus $((B \downarrow F),\mathtt{A})$ has a surinitial object in the fiber $(B \downarrow F)_{1_{\mathtt{A}(B)}}$. If we write down $B \xrightarrow{\eta_B} F(R_B)$ this surinitial object, then thanks to the proposition 9, it is initial in $(B \downarrow F)$. Then $F$ has a left adjoint: $G \dashv F$, and it is clearly a suradjoint. The converse is trivial. ■

## 7.4 A Theorem of Barr and Wells in the Surcategorical Context

As we are going to see,we have a surcategorical versus of the result that we can find in Borceux [6, proposition 4.3.6 page 206]. This theorem is a surcategorical adaptation of some results as we can find in Barr and Wells [2, § 9.3]. This theorem enables to conclude with a kind of cocompleteness of $T$-$\mathbb{C}$at and $CT$-$\infty$-$\mathbb{G}$r, which will show that some fusion hypotheses of our two suradjunctions are well-realized (see section 7.6).

**Theorem 3 (Barr-Wells's Surcategorical Theorem)** *Let* $(\mathscr{C},\mathtt{A})$ *a surcomplet and surcocomplet surcategory with* $K$*-equalizers. Let* $T$ *a surmonad on* $(\mathscr{C}, \mathtt{A})$*, which preserves* $\kappa$*-filtered surcolimits for some regular cardinal* $\kappa$*. In this case the surcategory* $(\mathscr{C}^T, \mathtt{A})$ *of suralgebras is surcomplet, surcocomplet, and has* $K$*-equalizers.* □

PROOF By using proposition 5, proposition 6, proposition 7 and theorem 2, there is no difficulty to transcript the proof from the classical case (as in [6]) to this surcategorical context. ■

## 7.5 Beck's Theorem in the Surcategorical Context

It is easy to see that Beck's theorem remains true in $\mathbb{C}at/\mathbb{G}$. This Beck's theorem enables to prove that the two important forgetful surfunctors ($U$ of theorem 10 and $V$ of theorem 11) are monadic (and more precisely surmonadic). We call this theorem "Sur-Beck's theorem" to refer to its surcategorical nature. Like in the classical case, we use two lemma which facilitate the demonstration of Sur-Beck's theorem (see [15]). But the proof of these two lemma and of Sur-Beck's theorem are very similar to the classical one (see [15]), and thus it is not necessary to give the details of the demonstrations. Contrary to the Freyd's suradjoint theorem and the Barr-Wells's surcategorical Theorem, we can notice that we do not need the presence of $K$-equalizers.

**Lemma 2** *Let* $(\mathscr{A}, \mathtt{A}) \underset{F}{\overset{G}{\underset{\longleftarrow}{\overset{\longrightarrow}{\top}}}} (\mathscr{X}, \mathtt{A})$, $(\mathscr{A}', \mathtt{A}) \underset{F'}{\overset{G'}{\underset{\longleftarrow}{\overset{\longrightarrow}{\top}}}} (\mathscr{X}, \mathtt{A})$*, two suradjunctions which have the same surmonad* $T$*. If we suppose that* $G$ *satisfied the hypothesis* 3 *of the theorem 4 then there is a unique surfunctor* $(\mathscr{A}', \mathtt{A}) \xrightarrow{M} (\mathscr{A}, \mathtt{A})$ *such as the following diagram commutes*

$$\begin{array}{ccccc}
\mathscr{X} & \xrightarrow{F'} & \mathscr{A}' & \xrightarrow{G'} & \mathscr{X} \\
{\scriptstyle 1_{\mathscr{X}}}\downarrow & & \downarrow{\scriptstyle M} & & \downarrow{\scriptstyle 1_{\mathscr{X}}} \\
\mathscr{X} & \xrightarrow[F]{} & \mathscr{A} & \xrightarrow[G]{} & \mathscr{X}
\end{array}$$

□

**Lemma 3** *In the situation* $(\mathscr{A},\mathtt{A}) \underset{F}{\overset{G}{\underset{\top}{\rightleftarrows}}} (\mathscr{X}\,\mathtt{A})$ *,* $G^T$ *creates surcoequalizers of* $(\mathscr{X}^T,\mathtt{A})$ *for absolute surcoequalizers, i.e given the diagram*

$$(x,h) \underset{d_1}{\overset{d_0}{\rightrightarrows}} (y,k)$$

*in one fiber of* $(\mathscr{X}^T,\mathtt{A})$ *such as the pair*

$$G^T((x,h)) \underset{G^T(d_1)}{\overset{G^T(d_0)}{\rightrightarrows}} G^T((y,k))\ ,\ (i.e\ \ x \underset{d_1}{\overset{d_0}{\rightrightarrows}} y\ )$$

*has an absolute surcoequalizer* $y \xrightarrow{e} z$*, so there is a unique* $T$*-algebra* $(z,m)$ *and a unique morphism* $(y,k) \xrightarrow{f} (z,m)$ *of* $(\mathscr{X}^T,\mathtt{A})$ *such as* $G^T(f) = e$ *and furthemore* $(y,k) \xrightarrow{f} (z,m)$ *is a surcoequalizer of the pair* $(x,h) \underset{d_1}{\overset{d_0}{\rightrightarrows}} (y,k)$ □

**Theorem 4 (Beck's Surcategorical Theorem)** *Let us Consider the suradjunction* $(\mathscr{A},\mathtt{A}) \underset{F}{\overset{G}{\underset{\top}{\rightleftarrows}}} (\mathscr{X},\mathtt{A})$ *with surmonad* $T$*, the canonical final suradjunction* $(\mathscr{X}^T,\mathtt{A}) \underset{F^T}{\overset{G^T}{\underset{\top}{\rightleftarrows}}} (\mathscr{X},\mathtt{A})$*, and the comparaison surfunctor* $(\mathscr{A},\mathtt{A}) \xrightarrow{K} (\mathscr{X}^T,\mathtt{A})$ *which is the unique surfunctor such as the following diagram commutes*

$$\begin{array}{ccccc}
\mathscr{X} & \xrightarrow{F} & \mathscr{A} & \xrightarrow{G} & \mathscr{X} \\
{\scriptstyle 1_{\mathscr{X}}}\downarrow & & \downarrow{\scriptstyle K} & & \downarrow{\scriptstyle 1_{\mathscr{X}}} \\
\mathscr{X} & \xrightarrow[F^T]{} & \mathscr{X}^T & \xrightarrow[G^T]{} & \mathscr{X}
\end{array}$$

*In this case the following conditions are equivalent*

1. *$K$ is an isomorphism in $\mathbb{C}at/\mathbb{G}$ (i.e there is a surfunctor $(\mathscr{X}^T,\mathtt{A})\overset{L}{\to}(\mathscr{A},\mathtt{A})$ such as $KL=1_{\mathscr{X}^T}$ and $LK=1_{\mathscr{A}}$).*

2. *$(\mathscr{A},\mathtt{A})\overset{G}{\to}(\mathscr{X},\mathtt{A})$ creates surcoequalizers of $a\underset{g}{\overset{f}{\rightrightarrows}}b$ for which the pair $G(a)\underset{G(g)}{\overset{G(f)}{\rightrightarrows}}G(b)$ has an absolute surcoequalizer.*

3. *$(\mathscr{A},\mathtt{A})\overset{G}{\to}(\mathscr{X},\mathtt{A})$ creates surcoequalizers of $a\underset{g}{\overset{f}{\rightrightarrows}}b$ for which the pair $G(a)\underset{G(g)}{\overset{G(f)}{\rightrightarrows}}G(b)$ has split surcoequalizers.* □

## 7.6 Fusion of Adjunctions

In section 9 and section 10 we build the two most important surmonads of this paper, the free surmonoids surmonad and the contractibility surmonad, and we need to do their "fusion" to obtain a new surmonad, which inherits at the same time properties of these two surmonads. This surmonad is the contractible surmonoids surmonad $\mathbb{B}=(B,\rho,b)$ of the theorem 1 which permits us to build the operads of non-strict cells. The fusion between suradjunctions require some hypotheses (see theorem 6) and naturally we shall see that our two suradjunctions fill these hypotheses.

The following "fusion theorem" is a generalization of techniques used by Batanin in [3]. This theorem is going to be shown especially powerful because the required hypotheses are very simple. As a result the fusion product of two adjunctions (which is a particular case of the fusion product of two suradjunctions) is possible under conditions that we can often run into. However it is the theorem 6, of which the fusion theorem is a consequence, that will be used to show theorem 1.

Here is a lemma, which is a particular case of the lemma 5:

**Lemma 4** *Let us consider the adjunction* $\mathscr{C} \underset{F}{\overset{U}{\underset{\top}{\rightleftarrows}}} \mathscr{B}$ *such as* $\mathscr{C}$ *has a co-equalizer and* $U$ *is faithful. Let the diagram* $B \overset{d_0}{\underset{d_1}{\rightrightarrows}} U(C)$ *in* $\mathscr{B}$*, then there is a unique morphism* $C \xrightarrow{q} Q$ *of* $\mathscr{C}$ *verifying* $U(q)d_0 = U(q)d_1$ *and which is universal for this property, i.e if we give ourselves another morphism* $C \xrightarrow{q'} Q'$ *of* $\mathscr{C}$ *such as* $U(q')d_0 = U(q')d_1$*, then there is a unique morphism* $Q \xrightarrow{h} Q'$ *of* $\mathscr{C}$ *such as* $U(h)U(q) = U(q')$. □

PROOF It is a particular case of lemma 5 where $\mathbb{G} = 1$ and where $1$ is the ponctual category. ■

The adjunction $\mathscr{C} \underset{F}{\overset{U}{\underset{\top}{\rightleftarrows}}} \mathscr{B}$ is fusionnable if over the ponctual category it is fusionnable too (See below the definition of fusionnable suradjunction over a category $\mathbb{G}$). Let us go to the fusion theorem.

**Theorem 5** *Let us consider the adjunction* $\mathscr{C} \underset{M}{\overset{U}{\underset{\top}{\rightleftarrows}}} \mathscr{B}$ *with monad* $(L, \mathfrak{m}, l)$*, and the adjunction* $\mathscr{D} \underset{H}{\overset{V}{\underset{\top}{\rightleftarrows}}} \mathscr{B}$ *with monad* $(C, m, c)$*. We suppose that these adjunctions are fusionnable. In this case, if we consider the cartesian square of categories*

$$\begin{array}{ccc} \mathscr{C} \times_{\mathscr{B}} \mathscr{D} & \xrightarrow{p_2} & \mathscr{D} \\ {\scriptstyle p_1}\downarrow & & \downarrow{\scriptstyle V} \\ \mathscr{C} & \xrightarrow[U]{} & \mathscr{B} \end{array}$$

*then the forgetful functor* $\mathscr{C} \times_{\mathscr{B}} \mathscr{D} \xrightarrow{O} \mathscr{B}$ *has a left adjoint:* $F \dashv O$. □

PROOF It is a particular case of the theorem 6 below where $\mathbb{G} = 1$ is the ponctual category. ■

**Lemma 5** *Consider the suradjunction* $(\mathscr{C},\mathtt{A}) \underset{F}{\overset{U}{\overset{\longrightarrow}{\underset{\longleftarrow}{\top}}}} (\mathscr{B},\mathtt{A})$ *such as* $\mathscr{C}$ *has surcoequalizers, and* $U$ *is faithful. Let the diagram* $B \underset{d_1}{\overset{d_0}{\rightrightarrows}} U(C)$ *in one fiber* $\mathscr{B}_G$ *above* $G \in \mathbb{G}$*, then there is a morphism* $C \xrightarrow{q} Q$ *of* $\mathscr{C}_G$ *verifying* $U(q)d_0 = U(q)d_1$ *and which is universal for this property, i.e if we give ourselves another morphism* $C \xrightarrow{q'} Q'$ *of* $\mathscr{C}$ *such as* $U(q')d_0 = U(q')d_1$*, then there is a unique morphism* $Q \xrightarrow{h} Q'$ *of* $\mathscr{C}$ *such as* $U(h)U(q) = U(q')$. □

PROOF This is similar to the proof in [3]. ■

Let the following suradjunction be: $(\mathscr{C},\mathtt{A}) \underset{F}{\overset{U}{\overset{\longrightarrow}{\underset{\longleftarrow}{\top}}}} (\mathscr{B},\mathtt{A})$ . It is fusionnable if the following properties are verified:

- $\mathscr{C}$ has surcoequalizers and $\overrightarrow{\mathbb{N}}$-surcolimits.
- $\mathscr{B}$ have $\overrightarrow{\mathbb{N}}$-surcolimits.
- $U$ is faithful and preserves $\overrightarrow{\mathbb{N}}$-surcolimits.

The theorem below is a key result for theorem 1:

**Theorem 6 (Surcategorical Fusion Theorem)** *Consider the suradjunction* $(\mathscr{C},\mathtt{A}) \underset{M}{\overset{U}{\overset{\longrightarrow}{\underset{\longleftarrow}{\top}}}} (\mathscr{B},\mathtt{A})$ *of underlying surmonad* $(L,\mathtt{m},l)$*, and the suradjunction* $(\mathscr{D},\mathtt{A}) \underset{H}{\overset{V}{\overset{\longrightarrow}{\underset{\longleftarrow}{\top}}}} (\mathscr{B},\mathtt{A})$ *of underlying surmonad* $(C,m,c)$*, such as they are both over the same category basis* $\mathbb{G}$ *and both are fusionnable. In this case, if we consider the cartesian square of categories*

$$\begin{array}{ccc} \mathscr{C} \times_{\mathscr{B}} \mathscr{D} & \xrightarrow{p_2} & \mathscr{D} \\ {\scriptstyle p_1}\downarrow & & \downarrow{\scriptstyle V} \\ \mathscr{C} & \xrightarrow[U]{} & \mathscr{B} \end{array}$$

*then the forgetful functor $\mathscr{C} \times_{\mathscr{B}} \mathscr{D} \xrightarrow{O} \mathscr{B}$ has a left adjoint $F \dashv O$. Besides this adjunction is a suradjunction. The surmonad of this suradjunction will be written $(B, \rho, b)$.* □

PROOF In this demonstration we often neglect the word "sur" because no confusion is possible.

Let $X \in \mathscr{B}(0)$ such as $\mathtt{A}(X) = G$. At first, we are going to build by induction an object $B(X)$ of $\mathscr{B}$ and secondly we shall reveal that $B(X)$ has got the expected universal property. We shall trivially see that $B(X)$ is in the fiber $\mathscr{B}_G$ and is built in it.

- If $n = 0$ we give ourselves the following diagram of $\mathscr{B}$:

$$C_0 = X \xrightarrow{l_0 = l(C_0)} L(C_0) \xrightarrow{\phi_0 = 1} L_0 \xrightarrow{c_0 = c(L_0)} C(L_0) \xrightarrow{\psi_0 = 1} C_1 \xrightarrow{l_1 = l(C_1)} L(C_1)$$

Thanks to the lemma, we obtain the morphism with the diagram $\phi_1$

$$L(C_0) \underset{d_1 = L(j_0) = L(\psi_0 c_0 \phi_0 l_0)}{\overset{d_0 = l_1 \psi_0 c_0 \phi_0}{\rightrightarrows}} L(C_1) \xrightarrow{\phi_1} L_1$$

What allows to extend the previous diagram

$$C_1 \xrightarrow{l_1} L(C_1) \xrightarrow{\phi_1} L_1 \xrightarrow{c_1 = c(L_1)} C(L_1)$$

And it allows again to obtain the morphism $\psi_1$

$$C(L_0) \underset{\delta_1 = C(k_0) = C(\phi_1 l_1 \psi_0 c_0)}{\overset{\delta_0 = c_1 \phi_1 l_1 \psi_0}{\rightrightarrows}} C(L_1) \xrightarrow{\psi_1} C_2$$

and thus to prolong once more the previous diagram

$$C_1 \xrightarrow{l_1} L(C_1) \xrightarrow{\phi_1} L_1 \xrightarrow{c_1} C(L_1) \xrightarrow{\psi_1} C_2 \xrightarrow{l_2} L(C_2)$$

We do an induction. We can suppose that up to the rank $n$ we can build these diagrams. In particular we give ourselves the following diagram

$$C_n \xrightarrow{l_n} L(C_n) \xrightarrow{\phi_n} L_n \xrightarrow{c_n} C(L_n) \xrightarrow{\psi_n} C_{n+1} \xrightarrow{l_{n+1}} L(C_{n+1})$$

where we especially note $j_n = \psi_n c_n \phi_n l_n$. We are going to show that we can prolong this type of diagram in the rank $n+1$. Thanks to the Lemma, we consider the morphism $\phi_{n+1}$

$$L(C_n) \underset{d_1=L(j_n)=L(\psi_n c_n \phi_n l_n)}{\overset{d_0=l_{n+1}\psi_n c_n \phi_n}{\rightrightarrows}} L(C_{n+1}) \xrightarrow{\phi_{n+1}} L_{n+1}$$

what allows to prolong the previous diagram

$$C_{n+1} \xrightarrow{l_{n+1}} L(C_{n+1}) \xrightarrow{\phi_{n+1}} L_{n+1} \xrightarrow{c_{n+1}=c(L_{n+1})} C(L_{n+1})$$

where we can particularly note $k_n = \phi_{n+1} l_{n+1} \psi_n c_n$. Then we consider, due to to the lemma, the morphism $\psi_{n+1}$

$$C(L_n) \underset{\delta_1=C(k_n)=C(\phi_{n+1} l_{n+1}\psi_n c_n)}{\overset{\delta_0=c_{n+1}\phi_{n+1} l_{n+1}\psi_n}{\rightrightarrows}} C(L_{n+1}) \xrightarrow{\psi_{n+1}} C_{n+2}$$

and thus to prolong still the previous diagram

$$C_{n+1} \xrightarrow{l_{n+1}} L(C_{n+1}) \xrightarrow{\phi_{n+1}} L_{n+1} \xrightarrow{c_{n+1}} C(L_{n+1}) \xrightarrow{\psi_{n+1}} C_{n+2} \xrightarrow{l_{n+2}} L(C_{n+2})$$

Thus for all $n \in \mathbb{N}$ we have this construction, what brings to light the filtered diagram built with these diagrams. This filtered diagram is noted $B_*$ and lives in the fiber $\mathscr{B}_G$. In particular the diagrams

$$\begin{array}{ccccc} L(C_{n-1}) & \underset{d_1=L(\psi_{n-1}c_{n-1}\phi_{n-1}l_{n-1})}{\overset{d_0=l_n\psi_{n-1}c_{n-1}\phi_{n-1}}{\rightrightarrows}} & L(C_n) & \underset{d_1=L(\psi_n c_n \phi_n l_n)}{\overset{d_0=l_{n+1}\psi_n c_n \phi_n}{\rightrightarrows}} & L(C_{n+1}) \\ & & \downarrow{\phi_n} & & \downarrow{\phi_{n+1}} \\ & & L_n & \xrightarrow[\lambda_n]{\cdots\cdots} & L_{n+1} \end{array}$$

show that

$\phi_{n+1}l_{n+1}\psi_n c_n\phi_n l_n\psi_{n-1}c_{n-1}\phi_{n-1} = \phi_{n+1}l_{n+1}\psi_n c_n\phi_n L(\psi_{n-1}c_{n-1}\phi_{n-1}$
$l_{n-1})$.

Thus according to the lemma, there is a unique morphism $L_n \xrightarrow{\lambda_n} L_{n+1}$, which is the forgetting of a morphism of $\mathscr{C}$, returning commutative these diagrams. Thus we obtain the filtered diagram $L_*$ of $\mathscr{B}$ which is the forgetting of a diagram filtered of $\mathscr{C}$ (and more precisely, is the forgetting of a diagram filtered of $\mathscr{C}_G$)

$$L_0 \xrightarrow{\lambda_0} L_1 \xrightarrow{\lambda_1} \cdots\cdots \to L_n \xrightarrow{\lambda_n} L_{n+1} \xrightarrow{\lambda_{n+1}} \cdots\cdots$$

where $B_*$ is an expanded diagram of $L_*$ i.e we have

$$\overbrace{C_0 \xrightarrow{l_0} L(C_0) \xrightarrow{\phi_0} L_*}^{B_*}$$

We also have the diagram

$$\begin{array}{ccccc} C(L_{n-2}) & \overset{\delta_0=c_{n-1}\phi_{n-1}l_{n-1}\psi_{n-2}}{\underset{\delta_1=C(\phi_{n-1}l_{n-1}\psi_{n-2}c_{n-2})}{\rightrightarrows}} & C(L_{n-1}) & \overset{\delta_0=c_n\phi_n l_n\psi_{n-1}}{\underset{\delta_1=C(\phi_n l_n\psi_{n-1}c_{n-1})}{\rightrightarrows}} & C(L_n) \\ & & \downarrow{\scriptstyle \psi_{n-1}} & & \downarrow{\scriptstyle \psi_n} \\ & & C_n & \overset{\kappa_n}{\cdots\cdots\to} & C_{n+1} \end{array}$$

which shows that

$\psi_n c_n\phi_n l_n\psi_{n-1}c_{n-1}\phi_{n-1}l_{n-1}\psi_{n-2} = \psi_n c_n\phi_n l_n\psi_{n-1}C(\phi_{n-1}l_{n-1}\psi_{n-2}$
$c_{n-2})$.

Thus according to the lemma, there is a unique morphism $C_n \xrightarrow{\kappa_n} C_{n+1}$ which is the forgetting of a morphism of $\mathscr{D}$ returning commutative these diagrams. Therefore we obtain the filtered diagram $C_*$ of $\mathscr{B}$ which is the forgetting of a diagram filtered of $\mathscr{D}$ (and more precisely, is the forgetting of a filtered diagram of $\mathscr{D}_G$)

$$C_1 \xrightarrow{\kappa_1} C_2 \xrightarrow{\kappa_2} \cdots\cdots > C_n \xrightarrow{\kappa_n} C_{n+1} \xrightarrow{\kappa_{n+1}} \cdots\cdots$$

where $B_*$ is an expanded diagram of $C_*$, i.e we have

$$\overbrace{C_0 \xrightarrow{c_0\phi_0 l_0} C(L_0) \xrightarrow{\psi_0} C_*}^{B_*}$$

Thus these diagrams $B_*$, $L_*$ and $C_*$ have the same surcolimit $B(X)$ in $\mathscr{B}$, and all the constructions are in the fiber $\mathscr{B}_G$. We put $L_* = U(M_*)$ and $M_* \to \Delta M_X$ its surcolimit surcocone (in $\mathscr{C}$), $C_* = V(H_*)$ and $H_* \to \Delta H_X$ its surcolimit surcocone (in $\mathscr{D}$). The surfunctors $U$ and $V$ preserving $\overrightarrow{\mathbb{N}}$-surcolimits, therefore $B(X)$ is the forgetting of the pair $(M_X, H_X)$ which is an object of $\mathscr{C} \times_{\mathscr{B}} \mathscr{D}$: $B(X) = O((M_X, H_X)) = U(M_X) = V(H_X)$. In particular $M_X$ is in the fiber $\mathscr{C}_G$ and $H_X$ is in the fiber $\mathscr{D}_G$. We put $F(X) = (M_X, H_X)$ which gives, as we are going to see, the desired left suradjoint of the forgetful surfunctor $O$, and where $(B, \rho, b)$ is the associated surmonad. $B(X)$ inherits at the same time the structure of the object $M_X$ (which lives in $\mathscr{C}$) and the structure of the object $H_X$ (which lives in $\mathscr{D}$). It is the reason why the surmonad $(B, \rho, b)$ can be called "fusion" of surmonads $(L, \mathfrak{m}, l)$ and $(C, m, c)$. We note $b_X$ the produced arrow $X \xrightarrow{b_X} B(X)$ The continuation consists in showing the universal character of $b_X$. We are going to show that if we give ourselves a morphism $X \xrightarrow{f} B_0$ of $\mathscr{B}$ such as $B_0$ is the forgetting of an object $(M_0, H_0)$ of $\mathscr{C} \times_{\mathscr{B}} \mathscr{D}$, then there is a unique morphism $(M_X, H_X) \xrightarrow{(h,k)} (M_0, H_0)$ of $\mathscr{C} \times_{\mathscr{B}} \mathscr{D}$ such as $O(h,k)b_X = f$. For that, we are going to use the filtered diagram $B_*$ with which we are going to build by induction a surcocone $B_* \to \Delta B_0$, and it will display the existence of the pair $(h,k)$.

– Let $g_0 = f$ and $f_0$ which is the extension of $f$ from $L_0 = L(X)$:

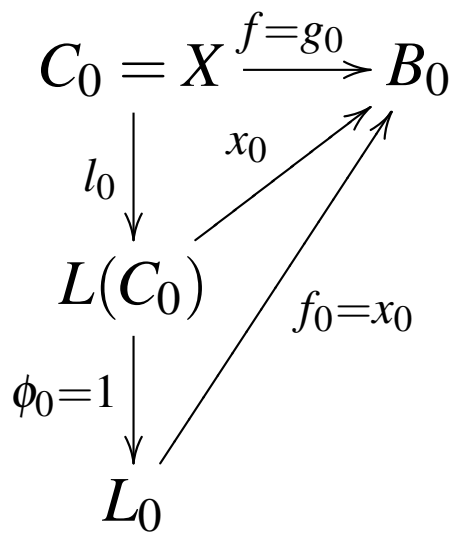

– We can suppose that this construction is up to the rank $n$. Thus in particular we have the following diagram

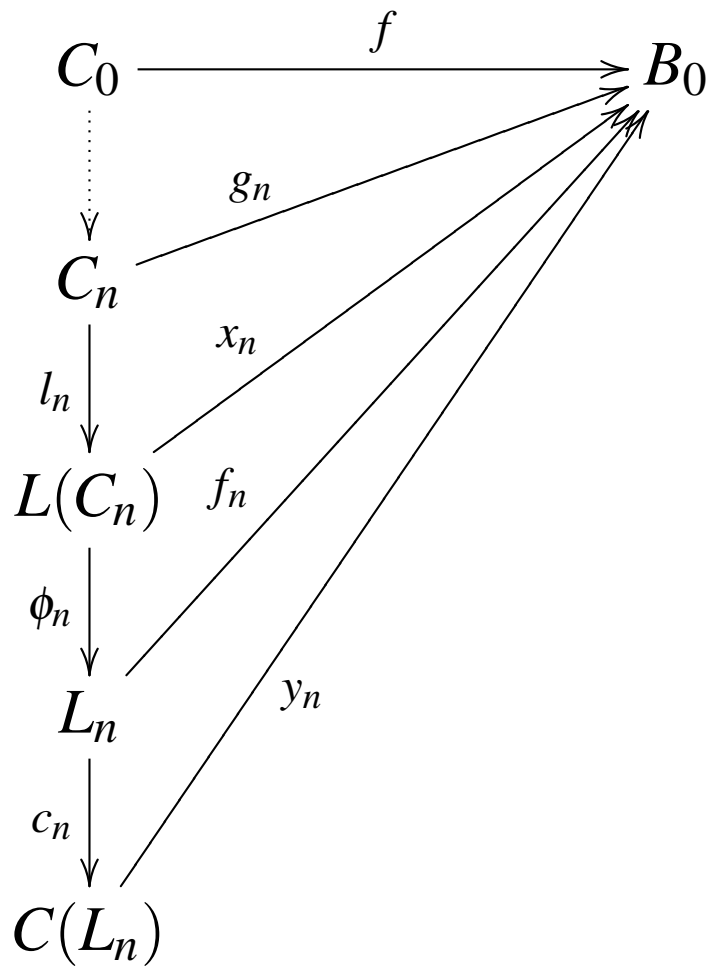

Also the natural transformation $1_{\mathscr{B}} \xrightarrow{c} C$ applied to

$$C(L_{n-1}) \xrightarrow{\phi_n l_n \psi_{n-1}} L_n$$

gives the equality

$$C(\phi_n l_n \psi_{n-1}) c(C(L_{n-1})) = c_n \phi_n l_n \psi_{n-1} = \delta_0$$

thus $y_n \delta_0 = y_n C(\phi_n l_n \psi_{n-1}) c(C(L_{n-1}))$. On the other hand

$$y_n \delta_0 = y_n \delta_0 m(L_{n-1}) c(C(L_{n-1}))$$

(unity axiom of monads), which leads to the equality

$$y_n C(\phi_n l_n \psi_{n-1}) = y_n \delta_0 m(L_{n-1})$$

(Do not forget that $y_n\delta_0$ is the forgetting of a morphism of $\mathscr{D}$ because $y_n\delta_0 = y_{n-1}$). What allows to write

$$\begin{aligned} y_n\delta_1 &= y_n C(k_{n-1}) = y_n C(\phi_n l_n \psi_{n-1} c_{n-1}) \\ &= y_n C(\phi_n l_n \psi_{n-1}) C(c_{n-1}) = y_n \delta_0 m(L_{n-1}) C(c_{n-1}) \\ &= y_n \delta_0 \text{ (unity axiom of monads)} \end{aligned}$$

So the universality of $\psi_n$ implies the existence of a unique morphism of $\mathscr{D}$ that the forgetting $g_{n+1}$ is such as $g_{n+1}\psi_n = y_n$. We also have the extension $x_{n+1}$ of $g_{n+1}$ from $L(C_{n+1})$. Then the natural transformation $1_{\mathscr{B}} \xrightarrow{l} L$ applied to $L(C_n) \xrightarrow{\psi_n c_n \phi_n} C_{n+1}$ gives the equality

$$L(\psi_n c_n \phi_n) l(L(C_n)) = l_{n+1} \psi_n c_n \phi_n = d_0$$

thus $x_{n+1}d_0 = x_{n+1}L(\psi_n c_n \phi_n) l(L(C_n))$, and

$$x_{n+1}d_0 = x_{n+1}d_0 \mathfrak{m}(C_n) l(L(C_n))$$

(unity axiom of monads), which leads to the equality

$$x_{n+1}L(\psi_n c_n \phi_n) = x_{n+1} d_0 \mathfrak{m}(C_n)$$

(do not forget that $x_{n+1}d_0$ is the forgetting of a morphism of $\mathscr{C}$ because $x_{n+1}d_0 = x_n$). What allows to write

$$\begin{aligned} x_{n+1}d_1 &= x_{n+1}L(j_n) = x_{n+1}L(\psi_n c_n \phi_n l_n) \\ &= x_{n+1}L(\psi_n c_n \phi_n) L(l_n) = x_{n+1} d_0 \mathfrak{m}(C_n) L(l_n) \\ &= x_{n+1} d_0 \text{ (unity axiom of monads)} \end{aligned}$$

Then the universality of $\phi_{n+1}$ implies the existence of a unique morphism of $\mathscr{C}$ which the forgetting $f_{n+1}$ is such as $f_{n+1}\phi_{n+1} = x_{n+1}$. We also have the extension $y_{n+1}$ of $f_{n+1}$ from $C(L_{n+1})$.

– Thus we obtain a cone $B_* \to \Delta B_0$, with $B_0 = O(M_0, H_0) = U(M_0) = V(H_0)$. We have the two surcocones as well $L_* \to \Delta U(M_0)$ and $C_* \to \Delta V(H_0)$. The functor $U$ preserving the $\overrightarrow{\mathbb{N}}$-surcolimits, the diagram of $\mathscr{B}$

$$\begin{array}{ccc} L_* & \longrightarrow & \Delta U(M_0) \\ & \searrow & \\ & & \Delta U(M_X) \end{array}$$

results of the diagram of $\mathscr{C}$

$$\begin{array}{ccc} M_* & \longrightarrow & \Delta M_0 \\ & \searrow & \\ & & \Delta M_X \end{array}$$

such as $M_* \to \Delta M_X$ is a surcolimit surcocone. There is consequently a unique morphism $h$ of $\mathscr{C}$ such as the triangle commutes

$$\begin{array}{ccc} M_* & \longrightarrow & \Delta M_0 \\ & \searrow & \uparrow{\scriptstyle \Delta h} \\ & & \Delta M_X \end{array} .$$

In the same way the surfunctor $V$ preserves $\overrightarrow{\mathbb{N}}$-surcolimits, so the diagram of $\mathscr{B}$

$$\begin{array}{ccc} C_* & \longrightarrow & \Delta V(H_0) \\ & \searrow & \\ & & \Delta V(H_X) \end{array}$$

results of the diagram of $\mathscr{D}$

$$\begin{array}{ccc} H_* & \longrightarrow & \Delta H_0 \\ & \searrow & \\ & & \Delta H_X \end{array}$$

such as $H_* \to \Delta H_X$ is a surcolimit surcocone. Therefore there is a unique morphism $k$ of $\mathscr{D}$ such as the following triangle commutes

$$\begin{array}{ccc} H_* & \longrightarrow & \Delta H_0 \\ & \searrow & \uparrow \vdots \, \Delta k \\ & & \Delta H_X \end{array} .$$

It shows the existence of the unique morphism $(h,k)$ of $\mathscr{C} \times_{\mathscr{B}} \mathscr{D}$ such as

$$\begin{array}{ccc} B_* & \longrightarrow & \Delta B_0 \\ & \searrow & \uparrow \vdots \, O(h,k) \\ & & \Delta B(X) \end{array} .$$

In consequence we obtain the morphism $(h,k)$ of $\mathscr{C} \times_{\mathscr{B}} \mathscr{D}$ such as $O(h,k)b_X = f$. Let $(h',k')$ another morphism of $\mathscr{C} \times_{\mathscr{B}} \mathscr{D}$ making the following triangle commute

$$\begin{array}{ccc} X & \xrightarrow{\quad f \quad} & B_0 = O(M_0, H_0) \\ {\scriptstyle b_X} \downarrow & \nearrow {\scriptstyle O(h',k')} & \\ B(X) = O(M_X, H_X) & & \end{array} .$$

We easily prove by induction that it makes commutative the following triangle of natural transformations

$$\begin{array}{ccc} B_* & \longrightarrow & \Delta B_0 \\ & \searrow & \uparrow \vdots \, O(h',k') \\ & & \Delta B(X) \end{array} .$$

and it immediatly prove the unicity of $(h,k)$.

Finally we obtain the following fusion diagram

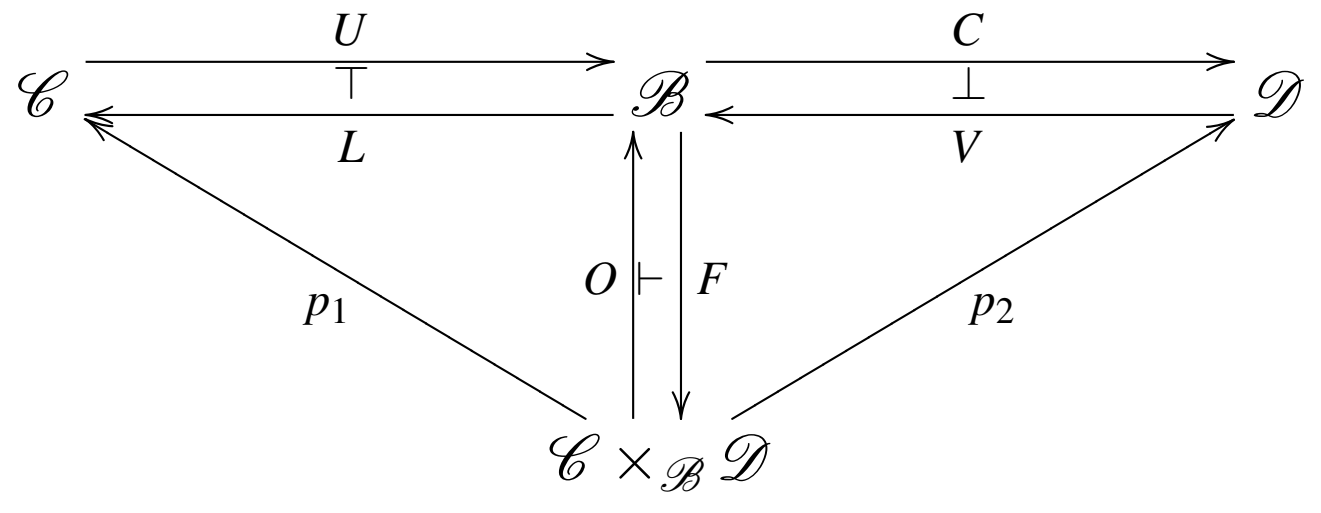

∎

# 8 Liberals Surmonoidals Surcategories

In [7] the authors suggest two constructions of the free monoid associated with an object of a monoidal category. This first construction (Bourn and Penon [7, proposition 1.2 page 14]) requires further properties on the underlying monoidal category that the authors call "numérale" for the surcategorical context (see Bourn and Penon [7, proposition 1.3.3 page 24]). The second construction of the free monoid such as it is found in Bourn and Penon [7, proposition 1.3 page 16] fits well with the pointed case (see section 8.2) and we are especially interested in this case (but in the surcategorical context). As the first construction, this second construction requires further properties on the underlying monoidal category. Therefore we call "liberal" those useful properties by which the free monoid can be obtained from this second construction. We shall make a small reminder of the main results but the reader is deeply encouraged to see the details of these constructions in [7] because we greatly use them at the end of the proof of the theorem 7. After we will show that all of these constructions apply in the surmonoidal context (which is the surcategorical versus of the monoidal context), where surmonoidals surcategories are for monoidals categories what surcategories are for categories. Although technics used here are close from those we find in [7], some concepts like Liberal Surmonoidal Surcategories and Pointed Surmonoidal

Surcategories are new. In particular the proof of the theorem 7 is similar to the proposition 10 below which is in [7].

## 8.1 Liberal Monoidal Categories

Let $\mathscr{V} = (\mathbb{V}, \bigotimes, I, u_l, u_r, aso)$ be a monoidal category. We sometimes design it by its underlyind category $\mathbb{V}$. $\mathscr{V}$ is liberal if the following properties hold:

- $\mathbb{V}$ has $\overrightarrow{\mathbb{N}}-$colimits and coequalizer;
- $\forall X \in \mathbb{V}(0)$, $(-)\bigotimes X$ and $X \bigotimes (-)$ preserves $\overrightarrow{\mathbb{N}}-$colimits;
- $\forall X \in \mathbb{V}(0)$, $(-)\bigotimes X$ preserves coequalizers.

Let $\mathbb{M}on(\mathbb{V})$ be the category of monoids associated with $\mathbb{V}$. We have a forgetful functor $\mathbb{M}on(\mathbb{V}) \xrightarrow{U} \mathbb{V}$, $(M, e, m) \longmapsto M$ and we have in Bourn and Penon [7, proposition 1.3 page 16]:

**Proposition 10** *If $\mathscr{V}$ is liberal and if $I$ is an initial object then the preceding forgetful functor has a left adjoint*

$$\mathbb{M}on(\mathbb{V}) \underset{Mo}{\overset{U}{\underset{\longleftarrow}{\overset{\longrightarrow}{\top}}}} \mathbb{V}$$

□

In order to construct this free monoid functor $\text{Mo}(-)$, we use the notion of graded monoid (defined in Bourn and Penon [7, page 12]). A graded monoid in a monoidal category $\mathscr{V}$ is given by a triple $((X_n)_{n\in\mathbb{N}}, (\iota_n)_{n\in\mathbb{N}}, (k_{n,m})_{n,m\in\mathbb{N}})$ where $(X_n)_{n\in\mathbb{N}}$ is a family of objects of $\mathscr{V}$, $(X_n \xrightarrow{\iota_n} X_{n+1})_{n\in\mathbb{N}}$ is a family of morphisms of $\mathscr{V}$, and $(X_n \bigotimes X_n \xrightarrow{k_{n,m}} X_{n+m})_{n,m\in\mathbb{N}}$ is a family of morphisms of $\mathscr{V}$, verifying some axioms that we can find in Bourn and Penon [7, page 12]. In [7] it is proved that every monoid has an underlying graded monoid and every graded monoid $((X_n)_{n\in\mathbb{N}}, (\iota_n)_{n\in\mathbb{N}}, (k_{n,m})_{n,m\in\mathbb{N}})$ is linked with a free monoid. Then the strategy to built the free monoid $\text{Mo}(X)$ for every $X \in \mathbb{V}$ is first to built a graded monoid $\Psi_X$ where this construction also requires the construction by induction of a secondary family of morphisms $(X \bigotimes X_n \xrightarrow{q_n} X_{n+1})_{n\in\mathbb{N}}$, then $\mathrm{Mo}(X)$ is also the free monoid associated with this graded monoid $\Psi_X = ((X_n)_{n\in\mathbb{N}}, (\iota_n)_{n\in\mathbb{N}}, (k_{n,m})_{n,m\in\mathbb{N}})$.

We must remember that $\mathrm{Mo}(X) = \mathrm{colim} X_n$ and $(X_n)_{n\in\mathbb{N}}$ is built by induction with morphisms $X_n \xrightarrow{\iota_n} X_{n+1}$, $X \bigotimes X_n \xrightarrow{q_n} X_{n+1}$, by considering the coequalizer $q_{n+1} := coker(y_n^0, y_n^1)$, where

- $y_n^0 = Id \otimes \iota_n$,
- $y_n^1 = (\ X \bigotimes X_n \xrightarrow{q_n} X_{n+1} \xrightarrow{u_l^{-1}} I \bigotimes X_{n+1} \xrightarrow{!\otimes Id} X \bigotimes X_{n+1}\ )$,
- $\iota_{n+1} = (\ X_{n+1} \xrightarrow{u_l^{-1}} I \bigotimes X_{n+1} \xrightarrow{!\otimes Id} X \bigotimes X_{n+1} \xrightarrow{q_{n+1}} X_{n+2}\ )$,

and where the initialization is given by $X_0 = I, X_1 = X$, $I \xrightarrow{\iota_0 = !_X} X$, $X \bigotimes I \xrightarrow{q_0 = u_r} X$.

Morphims $k_{n,m}$ are built by induction (see Bourn and Penon [7, page 16 and page 17]), but we do not describe it here because we do not explicitly need them anymore. Let $(X_n \xrightarrow{l_n} \mathrm{Mo}(X))_{n\in\mathbb{N}}$, the universal cocone defining $\mathrm{Mo}(X)$. The associated universal arrow is $X \xrightarrow{l(X)=l_1} \mathrm{Mo}(X)$. Let us remind that the multiplication $\mathrm{Mo}(X) \bigotimes \mathrm{Mo}(X) \xrightarrow{m} \mathrm{Mo}(X)$ is the unique arrow such as $\forall n,m \in \mathbb{N}$: $m(l_n \otimes l_m) = l_{n+m} k_{n,m}$. When $n = 1$ we have $k_{1,m} = q_m$ which gives the equality $m(l_1 \otimes l_m) = l_{m+1} q_m$ and which will be useful for the construction of the free surmonoid (see result 2).

## 8.2 Liberal Surmonoidal Surcategories

Let $\mathbb{G}$ be some fixed category.

We shall expand further on the "surmonoidal" context, what we have made for the monoidal context. In a way we are going to point out that the results of $\mathbb{C}$at, which has enabled to build the free contractible operad of non-strict $\infty$-categories of Batanin (see [3]) are true in $\mathbb{C}\mathrm{at}/\infty\text{-}\mathbb{G}\mathrm{r}$, which are going to give us the free contractibles colored operads of Batanin's non-strict cells.

Let us now briefly recall the definition of surmonoidal surcategory.

Let $\mathbb{G}$ be a fixed category. A surmonoidal surcategory (over $\mathbb{G}$) is a monoidal object of the 2-category $\mathbb{C}\mathrm{at}/\mathbb{G}$. A surmonoidal surcategory is thus given by a 7-uple: $\mathscr{E}=(\mathbb{E},\mathtt{A},\bigotimes,I,u_l,u_r,aso)$ where:

- $\mathtt{A}$ is a functor: $\mathbb{E}\xrightarrow{\mathtt{A}}\mathbb{G}$;
- $(\mathbb{E}\times_{\mathbb{G}}\mathbb{E},\mathtt{A})\xrightarrow{\bigotimes}(\mathbb{E},\mathtt{A})$ is a morphism of $\mathbb{C}\mathrm{at}/\mathbb{G}$, where $\mathbb{E}\times_{\mathbb{G}}\mathbb{E}$ results of the cartesian square given by $\mathtt{A}$;
- $\mathbb{G}\xrightarrow{I}\mathbb{E}$ is a functor and a section (i.e we have $\mathtt{A}I=1_{\mathbb{G}}$);
- $u_r$ and $u_l$ are natural isomorphisms: $\bigotimes(1_{\mathbb{E}},I\mathtt{A})\xrightarrow{u_r}1_{\mathbb{E}}$, $\bigotimes(I\mathtt{A},1_{\mathbb{E}})\xrightarrow{u_l}1_{\mathbb{E}}$;
- $aso$ is a natural isomorphism: $\bigotimes(\bigotimes\times 1_{\mathbb{E}})\xrightarrow{aso}\bigotimes(1_{\mathbb{E}}\times\bigotimes)$.

And these data satisfy the usual conditions of coherences i.e those given by the axioms of monoidal categories. A simple consequence of this definition is that for every object $B$ of $\mathbb{G}$ each fiber $\mathbb{E}_B$ is a monoidal category. We write with the same notation in each fiber the tensor product because the context will prevent any confusion.

**Remark 6** Obviously, strict surmonoidal surcategories are such as $u_l,u_r$ and $aso$ are natural identities. □

Let $\mathscr{E}=(\mathbb{E},\mathtt{A},\bigotimes,I,u_l,u_r,aso)$ and $\mathscr{E}'=(\mathbb{E}',\mathtt{A}',\bigotimes',I',u'_l,u'_r,aso')$ be two surmonoidal surcategories with respective base categories $\mathbb{G}$ and $\mathbb{G}'$. A strict morphism $\mathscr{E}\xrightarrow{(F,F_0)}\mathscr{E}'$, is given by two functors $\mathbb{E}\xrightarrow{F}\mathbb{E}'$ and $\mathbb{G}\xrightarrow{F_0}\mathbb{G}'$ such as $F_0\mathtt{A}=\mathtt{A}'F$, $FI=I'F_0$ and $F\bigotimes=\bigotimes'(F\times_{F_0}F)$.

Let $\mathscr{E}$ be a surmonoidal surcategory. A surmonoid in $\mathscr{E}$ is given by a pair $(\mathscr{C};C_0)$ where $C_0\in\mathbb{G}(0)$ and $\mathscr{C}=(C,m,e)$ is a monoid in $\mathbb{E}_{C_0}$ ($m$ is the multiplication and $e$ is the unity). Thus $(\mathscr{C};C_0)$ can be noted as well $(C,m,e;C_0)$.

If $(\mathscr{C};C_0)$ and $(\mathscr{C}';C'_0)$ are surmonoids, a morphism

$$(\mathscr{C};C_0) \xrightarrow{(f,f_0)} (\mathscr{C}';C'_0),$$

is given by a pair $(f,f_0)$ where $C_0 \xrightarrow{f_0} C'_0$ is an arrow in $\mathbb{G}$ and $\mathscr{C} \xrightarrow{f} \mathscr{C}'$ is given by an arrow $C \xrightarrow{f} C'$ in $\mathbb{E}$ such as $\mathtt{A}(f) = f_0$ and $fm = m'(f \otimes_{f_0} f)$, $fe = e'I(f_0)$. We note $/\mathbb{M}on(\mathbb{E},\mathtt{A})$ the category of surmonoids of $\mathscr{E}$.

Let $\mathscr{E}$ be a surmonoidal surcategory. It is liberal if the following two conditions are checked

- $\forall B \in \mathbb{G}(0)$, the fiber $\mathbb{E}_B$ is a liberal monoidal category.
- $\forall B \in \mathbb{G}(0)$, the canonical inclusion functor $\mathbb{E}_B \hookrightarrow \mathbb{E}$ preserves coequalizer and $\overrightarrow{\mathbb{N}}$-colimits.

Let $(\mathscr{C};C_0)$ be a surmonoid of $\mathscr{E}$, then

**Proposition 11** *The pair $(\mathbb{E}/C,\widehat{\mathtt{A}})$, such as $\mathbb{E}/C \xrightarrow{\widehat{\mathtt{A}}} \mathbb{G}/\mathtt{A}(C)$, $x \longmapsto \mathtt{A}(x)$, produces a surmonoidal surcategory*

$$\mathscr{E}/C = (\mathbb{E}/C,\widehat{\mathtt{A}},\widehat{\bigotimes},\widehat{I},\widehat{u_l},\widehat{u_r},\widehat{aso}) \qquad \square$$

The proof is in Bourn and Penon [7, page 22] but let us remind that if $(X,x),(Y,y) \in \mathbb{E}/C$ then $(X,x)\widehat{\otimes}(Y,y) := (X \bigotimes Y, m(x \otimes y))$. If $b \in \mathbb{G}/\mathtt{A}(C)$ then $\widehat{I}(b) := eI(b)$. The 2-cells $\widehat{u_l}$, $\widehat{u_r}$, $\widehat{aso}$ are also provided with the corresponding data of $\mathscr{E}$.

When surcoequalizers exist in $\mathscr{E}$, it is not difficult to see that surcoequalizers in $\mathbb{E}/C$ are computed by it, and we have the same phenomenon for $\overrightarrow{\mathbb{N}}$-surcolimits. So we have the following easy proposition that is left for the reader.

**Proposition 12** *If* $\mathscr{E} = (\mathbb{E}, \mathtt{A}, \bigotimes, I, u_l, u_r, aso)$ *is a liberal surmonoidal surcategory then* $\mathscr{E}/C = (\mathbb{E}/C, \widehat{\mathtt{A}}, \widehat{\bigotimes}, \widehat{I}, \widehat{u}_l, \widehat{u}_r, \widehat{aso})$ *is a liberal surmonoidal surcategory, and the morphism* $\mathscr{E}/C \xrightarrow{(S,S_0)} \mathscr{E}$ *given by the functor* $\mathbb{E}/C \xrightarrow{S} \mathbb{E}$, $(X,x) \longmapsto X$, *is a strict morphism of surmonoidal surcategories which preserves the liberal structure.* □

We have the following proposition too

**Proposition 13** *If* $\mathscr{E}$ *is a liberal surmonoidal surcategory and if* $(\mathscr{C};C_0) \xrightarrow{(h,h_0)} (\mathscr{C}';C'_0)$ *is a morphism of surmonoids, then the morphism*

$$\mathscr{E}/C \xrightarrow{(h^*,h_0^*)} \mathscr{E}/C'$$

*is a strict morphism of surmonoidal surcategories which preserves the liberal structure.* □

PROOF The fact that $(h^*, h_0^*)$ is a strict morphism of surmonoidal surcategories has already been shown (Bourn and Penon [7, page 25]) and the fact that $h^*$ preserves $\overrightarrow{\mathbb{N}}$-surcolimits has already been proved for the numeral context [see 7, page 25]. We only have to show that $h^*$ preserves surcoequalizers and it is evident by construction. ■

Now we have some materials to show the main theorem of this paragraph.

**Theorem 7** *Let* $\mathscr{E} = (\mathbb{E}, \mathtt{A}, \bigotimes, I, u_l, u_r, aso)$ *be a liberal surmonoidal surcategory such as* $\forall B \in \mathbb{G}(0)$ *the object* $I(B)$ *is initial in the fiber* $\mathbb{E}_B$ *and such as* $\forall b \in \mathbb{G}(1)$ *the object* $I(b)$ *is initial in the fiber* $\mathbb{E}_b$*, then the forgetful surfunctor*

$$(/\mathbb{M}on(\mathbb{E},\mathtt{A}),\mathtt{A}) \xrightarrow{U} (\mathbb{E},\mathtt{A})$$

*has a left suradjoint* $M \dashv U$ *and it is surmonadic.* □

PROOF It is similar to the proof of proposition 10 and we just need to adapt it in the surcategorical context. In particular we use proposition 10, the reminders in section 8.1, the proposition 13 plus the following two results (The first result below is a refinement of the proposition 13. We prove these two results by induction):

**Result 1** *Let* $(\mathscr{C};C_0) \xrightarrow{(h,h_0)} (\mathscr{C}';C'_0)$ *be a morphism of surmonoids, then if* $(X,x) \in \mathbb{E}/C$ *then* $\forall n \in \mathbb{N}$*,* $h^*((X,x)_n) = (h^*((X,x)))_n$*, where* $(X,x)_n$ *is the* $n^{th}$ *object of the graded monoid associated with* $(X,x)$ *(see [7, 1.2.3 page 12] for definition and results about graded monoid).*

**Result 2** $\forall n \in \mathbb{N}$*,* $(X,l(X))_n = (X_n,l_n)$*, where* $X_n \xrightarrow{l_n} \mathrm{Mo}(X)$ *is an arrow of the colimit cocone defining* $\mathrm{Mo}(X)$*, and where* $(X,l(X))_n$ *is the* $n^{th}$ *object of the graded monoid associated with* $(X,l(X))$.

The surmonadicity of $U$ is a simple consequence of theorem 4. In particular this surmonadicity has already been proved in the numeral context [7, see proposition 1.3.1, page 20]. ■

Now we can study the important case of pointed surmonoidal surcategories. In [7] it is proved that to any monoidal category $\mathscr{V} = (\mathbb{V},\bigotimes,I,u_l,u_r,aso)$ we associate its pointed monoidal category

$$Pt(\mathscr{V}) = (Pt(\mathbb{V}),\widetilde{\bigotimes},\widetilde{I},\widetilde{u_l},\widetilde{u_r},\widetilde{aso})$$

and if $\mathscr{V}$ was liberal then $Pt(\mathscr{V})$ remained liberal.

We can expand to the surmonoidal context this construction and this result. Let $\mathscr{E} = (\mathbb{E},\mathtt{A},\bigotimes,I,u_l,u_r,aso)$ be a surmonoidal surcategory over a fixed category $\mathbb{G}$.

Let $Pt(\mathbb{E})$ the category with objects the pairs $(X,x)$ where $X \in \mathbb{E}(0)$ and $I(\mathtt{A}(X)) \xrightarrow{x} X \in \mathbb{E}(1)$, and which has for arrows $(X,x) \xrightarrow{f} (Y,y)$, given by morphism $X \xrightarrow{f} Y$ of $\mathbb{E}$ such as $fx = yI_{\mathtt{A}(f)}$. In this case we have **Proposition 14** *The pair* $(Pt(\mathbb{E}),\widetilde{\mathtt{A}})$ *such as:* $(X,x) \overset{\widetilde{\mathtt{A}}}{\longmapsto} \mathtt{A}(X)$ *, produces a structure of surmonoidal surcategory*

$$Pt(\mathscr{E}) = (Pt(\mathbb{E}),\widetilde{\mathtt{A}},\widetilde{\bigotimes},\widetilde{I},\widetilde{u_l},\widetilde{u_r},\widetilde{aso}) \qquad \square$$

PROOF • Its tensor is the bifunctor $Pt(\mathbb{E}) \times_{\mathbb{G}} Pt(\mathbb{E}) \overset{\widetilde{\otimes}}{\longrightarrow} Pt(\mathbb{E})$, $((X,x),(Y,y)) \longmapsto (X,x)\widetilde{\otimes}(Y,y) := (X \bigotimes Y, (x\otimes y)u_l^{-1})$.

- Its "unity" functor is $\mathbb{G} \overset{\widetilde{I}}{\to} Pt(\mathbb{E})$, $G \longmapsto (I(G), 1_{I(G)})$.

- Left and right isomorphisms of unity: For all $(X,x)$ of $Pt(\mathbb{E})(0)$ the tensor $\widetilde{I}(\widetilde{\mathtt{A}}(X,x))\widetilde{\otimes}(X,x)$ is given by the morphism $(1_{I(\mathtt{A}(X))} \otimes x)u_l^{-1}$ of $\mathbb{E}$, and we have $u_l(1_{I(\mathtt{A}(X))} \otimes x)u_l^{-1} = x$ thanks to the equality

  $$u_l(X)(1_{I(\mathtt{A}(X))} \otimes x) = xu_l(I(\mathtt{A}(X))).$$

  Thus we get

  $$\widetilde{I}(\widetilde{\mathtt{A}}(X,x))\widetilde{\otimes}(X,x) \xrightarrow{\widetilde{u_l}(X,x)} (X,x)$$

  and $\widetilde{u_l}(X,x)$ given by $u_l(X)$ is a good candidate to define $\widetilde{u_l}$. Thus we obtain the natural transformation $\widetilde{\otimes}(\widetilde{I}\widetilde{\mathtt{A}},Id) \overset{\widetilde{u_l}}{\longrightarrow} Id$ which is in fact, an underlying datum of its 2-cell $\widetilde{u_l}$. In the same way we obtain the 2-cell $\widetilde{\otimes}(Id,\widetilde{I}\widetilde{\mathtt{A}}) \overset{\widetilde{u_r}}{\Longrightarrow} Id$ .

- The tensor products

  $$((X,x)\widetilde{\bigotimes}(Y,y))\widetilde{\bigotimes}(Z,z) \text{ and } (X,x)\widetilde{\bigotimes}((Y,y)\widetilde{\bigotimes}(Z,z))$$

  are respectively given by

  $$[((x\otimes y)u_l^{-1}) \otimes z]u_l^{-1} \text{ and } [x\otimes((y\otimes z)u_l^{-1})]u_l^{-1},$$

and we have the equality

$$aso[((x\otimes y)u_l^{-1})\otimes z]u_l^{-1} = [x\otimes((y\otimes z)u_l^{-1})]u_l^{-1}$$

due to the naturality of *aso* and the underlying surmonoid structure of $I(\mathtt{A}(X))$. We consequently obtain

$$((X,x)\widetilde{\bigotimes}(Y,y))\widetilde{\bigotimes}(Z,z) \xrightarrow{\widetilde{aso}} (X,x)\widetilde{\bigotimes}((Y,y)\widetilde{\bigotimes}(Z,z))$$

where in particular $\widetilde{aso}$ is given by *aso*, and is the good candidate to be the 2-cells of associativity. Thus we obtain the natural transformation $\widetilde{\bigotimes}(\widetilde{\bigotimes}\times Id) \xrightarrow{\widetilde{aso}} \widetilde{\bigotimes}(Id\times\widetilde{\bigotimes})$ which in reality is an underlying datum of the 2-cell $\widetilde{aso}$, and with this description of $Pt(\mathscr{E})$ it is now easy to see that it is a surmonoidal surcategory. ∎

As for $\mathscr{E}/C$, when surcoequalizers exist in $\mathscr{E}$, then we can see that surcoequalizers in $Pt(\mathscr{E})$ are computed by it, and we have the same phenomenon for $\overrightarrow{\mathbb{N}}$- surcolimits. So we have the following easy proposition intended for the readers.

**Proposition 15** *If* $\mathscr{E} = (\mathbb{E},\mathtt{A},\bigotimes,I,u_l,u_r,aso)$ *is a liberal surmonoidal surcategory then* $Pt(\mathscr{E}) = (Pt(\mathbb{E}),\widetilde{\mathtt{A}},\widetilde{\bigotimes},\widetilde{I},\widetilde{u}_l,\widetilde{u}_r,\widetilde{aso})$ *stays a liberal surmonoidal category, and trivially the functor* $\widetilde{I}$ *send objects and arrows of* $\mathbb{G}$ *to initials objects in the corresponding fibers.* □

The following proposition is easy and do not thus require detailed proof. It is the surmonoidal version of the result in Bourn and Penon [7, 1.2.1 page 10].

**Proposition 16** *If* $\mathscr{E} = (\mathbb{E},\mathtt{A},\bigotimes,I,u_l,u_r,aso)$ *is a surmonoidal surcategory, then we have the commutative triangle*

$$\begin{array}{ccc} /\mathbb{M}on(\mathbb{E},\mathtt{A}) & \xrightarrow{\varphi} & /\mathbb{M}on(Pt(\mathbb{E}),\widetilde{\mathtt{A}}) \\ & \searrow^{U} & \downarrow U' \\ & & Pt(\mathbb{E}) \end{array}$$

*such as $\varphi$ is an isomorphism given by $\varphi((C,e,m;C_0)) = ((C,e),e,m;C_0)$, and with $U((C,e,m;C_0)) = (C,e)$ and $U'(((C,x),e,m;C_0)) = (C,e)$.* □

With the theorem and the previous propositions we have at once:

**Theorem 8** *If $\mathscr{E} = (\mathbb{E},\mathtt{A},\bigotimes,I,u_l,u_r,aso)$ is a liberal surmonoidal surcategory then the forgetful surfunctor*

$$(/\mathbb{M}on(\mathbb{E},\mathtt{A}),\mathtt{A}) \xrightarrow{U} (Pt(\mathbb{E}),\widetilde{\mathtt{A}})\ ,\ (C,e,m;C_0) \longmapsto (C,e)$$

*has a left suradjoint and is surmonadic.* □

**Remark 7** Let us note $M$ the left suradjoint of $U$, then if we applied $(X,x) \in Pt(\mathbb{E})$ to the unity $1_{Pt(\mathbb{E})} \xrightarrow{l} UM$ of this suradjunction, we obtain the morphism $(X,x) \xrightarrow{l((X,x))} U(M(X,x))$ of $Pt(\mathbb{E})$ i.e $(X,x) \xrightarrow{l((X,x))} U(\overline{X},e,m;X_0) = (\overline{X},e)$. And in particular this morphism gives us the equality $l((X,x))x = e$. This equality is important because it shows, in the particular context of colored operads of the previous paragraphs, that the operads of non-strict cells are well-provided with a system of operations. □

# 9 Obtaining of the First Monad

With the notations of the paragraph section 1.1 we have

**Proposition 17** *The pair $(T\text{-}\mathbb{G}r,\mathtt{A})$*

$$T\text{-}\mathbb{G}r \xrightarrow{\mathtt{A}} \infty\text{-}\mathbb{G}r\ ,\ (C,d,c) \longmapsto \text{codomain}(c)$$

*has a structure of liberal surmonoidal surcategory.* □

If $\mathscr{C}$ is a topos we shall note $\mathscr{C}/B \xrightarrow{f^*} \mathscr{C}/A$ the pullback functor associated with an arrow $A \xrightarrow{f} B$ of $\mathscr{C}$, and $\mathscr{C}/A \xrightarrow{\Sigma_f} \mathscr{C}/B$ the composition functor. We have the usual adjunctions: $\Sigma_f \dashv f^* \dashv \pi_f$, where $\pi_f$ is the internal product functor.

PROOF The tensor $\otimes$ of $(T\text{-}\mathbb{G}\mathrm{r},\mathtt{A})$ is defined in the following way ($p_1$ and $p_2$ are the obvious projections)

$$T\text{-}\mathbb{G}\mathrm{r}\times_{\infty\text{-}\mathbb{G}\mathrm{r}}T\text{-}\mathbb{G}\mathrm{r}\xrightarrow{\otimes}T\text{-}\mathbb{G}\mathrm{r}\,,$$

$$((C,d,c),(C',d',c'))\longmapsto(T(C)\times_{T(G)}C',\mu(G)T(d)\pi_0,c'\pi_1)$$

where $G=\mathtt{A}((C,d,c))=\mathtt{A}((C',d',c'))$, and where in particular the tensorial product of morphisms

$$\Big((C,d,c)\xrightarrow{f}(D,d,c)\Big)\otimes\Big((C',d',c')\xrightarrow{f'}(D',d',c')\Big)$$

is obtained due to the cartesian square given by $(D,d,c)\otimes(D',d',c')$, i.e it is given by the morphism $C\otimes C'\xrightarrow{T(f)\times_{T(g)}f'}D\otimes D'$ where $g=\mathtt{A}(f)=\mathtt{A}(f')$. We also have the unity functor $\infty\text{-}\mathbb{G}\mathrm{r}\xrightarrow{I}T\text{-}\mathbb{G}\mathrm{r}$, given by: $G\longmapsto(G,l(G),1_G)$.

The natural isomorphism $\otimes\circ(1_{T\text{-}\mathbb{G}\mathrm{r}},I\mathtt{A})\xrightarrow{u_r}1_{T\text{-}\mathbb{G}\mathrm{r}}$ is trivial, and the natural isomorphism $\otimes\circ(I\mathtt{A},1_{T\text{-}\mathbb{G}\mathrm{r}})\xrightarrow{u_l}1_{T\text{-}\mathbb{G}\mathrm{r}}$ results from the cartesianity of $T$. We easily have the natural isomorphism as well

$$\bigotimes\circ(\bigotimes\times 1_{T\text{-}\mathbb{G}\mathrm{r}})\xrightarrow{aso}\bigotimes\circ(1_{T\text{-}\mathbb{G}\mathrm{r}}\times\bigotimes)$$

and we show without any problems that $(T\text{-}\mathbb{G}\mathrm{r},\mathtt{A})$ is a surmonoidal surcategory.

Let us now prove its liberal properties. Let $G\in\infty\text{-}\mathbb{G}\mathrm{r}$, and let $(C,d,c)\in T\text{-}\mathbb{G}\mathrm{r}_G\simeq\infty\text{-}\mathbb{G}\mathrm{r}/T(G)\times G$, $(D,d',c')\mapsto(D,(d',c'))$. In this case we show that we have the isomorphism of functors

$$(\Sigma_{(\mu(G)T(d)\times 1_G)})(T(c)\times 1_G)^*\simeq(-)\bigotimes(C,d,c).$$

And if we put $\widehat{T}(C,d,c):=(T(C),T(d),T(c))$ and if $(D,d',c')\in T\text{-}\mathbb{G}\mathrm{r}_G$ then we show that we have the isomorphism of functors

$$(\Sigma_{\mu(G)\times c'})(1_{T^2(G)}\times d')^*\widehat{T}\simeq(D,d',c')\bigotimes(-).$$

With these descriptions of functors $(-)\otimes(C,d,c)$ and $(D,d',c')\otimes(-)$, and knowing that the category $\infty\text{-}\mathbb{G}\mathrm{r}/T(G)\times G$ is a topos, the liberal structure is then easy:

- $(-)\bigotimes(C,d,c)$ preserves $\overrightarrow{\mathbb{N}}$−colimits and coequalizers.
- $(D,d',c')\bigotimes(-)$ preserves $\overrightarrow{\mathbb{N}}$−colimits in particular because $T$ trivially preserves them (and thus $\widehat{T}$ preserves them).
- We can see without any difficulties that the canonical inclusion of the fiber above $G$: $T\text{-}\mathbb{G}\mathrm{r}_G \hookrightarrow T\text{-}\mathbb{G}\mathrm{r}$, preserves $\overrightarrow{\mathbb{N}}$−colimits and coequalizers. ■

Thus thanks to the theorem 8

**Theorem 9** *The forgetful surfunctor*

$$(T\text{-}\mathbb{C}at,\mathtt{A}) \xrightarrow{U} (T\text{-}\mathbb{G}r_p,\mathtt{A}) ,$$

*where:* $(T\text{-}\mathbb{C}at,\mathtt{A}) := (/\mathbb{M}on(T\text{-}\mathbb{G}r,\mathtt{A}),\mathtt{A})$ *and* $(T\text{-}\mathbb{G}r_p,\mathtt{A}) := (Pt(T\text{-}\mathbb{G}r),\mathtt{A})$*, has a left suradjoint* $M \dashv U$*. Moreover this forgetful surfunctor* $U$ *is surmonadic.* □

We have $(T\text{-}\mathbb{G}\mathrm{r}_c,\mathtt{A}_c) = i^*(T\text{-}\mathbb{G}\mathrm{r},\mathtt{A})$ where $\infty\text{-}\mathbb{G}\mathrm{r}_c \overset{i}{\hookrightarrow} \infty\text{-}\mathbb{G}\mathrm{r}$ is the inclusion functor, and we know that the functor $i^*$ preserves many things (in particular the surmonoidal structure and the liberal structure), which also prove that we have the

**Proposition 18** *The pair* $(T\text{-}\mathbb{G}r_c,\mathtt{A}_c)$

$$T\text{-}\mathbb{G}r_c \xrightarrow{\mathtt{A}_c} \infty\text{-}\mathbb{G}r_c\ ,\ (C,d,c) \longmapsto \mathrm{codomain}(c)$$

*have a structure of liberal surmonoidal surcategory.* □

We can now establish the existence of the free surmonoids surmonad, the first important surmonad of this article.

**Theorem 10** *The forgetful surfunctor*

$$(T\text{-}\mathbb{C}at_c,\mathtt{A}) \xrightarrow{U} (T\text{-}\mathbb{G}r_{c,p},\mathtt{A}) ,$$

*where:*

$$(T\text{-}\mathbb{C}at_c, \mathtt{A}) := (/\mathbb{M}on(T\text{-}\mathbb{G}r_c, \mathtt{A}_c), \mathtt{A})$$

$$\text{and } (T\text{-}\mathbb{G}r_{c,p}, \mathtt{A}) := (Pt(T\text{-}\mathbb{G}r_c), \mathtt{A}),$$

*has a left suradjoint $M \dashv U$. Besides this forgetful surfunctor $U$ is surmonadic.* □

PROOF It is an application of the theorem 8 and the previous proposition 18.■

# 10 Obtaining of the Second Monad

Categories $T$-$\mathbb{G}\mathrm{r}_c$, $CT$-$\mathbb{G}\mathrm{r}_c$, $T$-$\mathbb{G}\mathrm{r}_{c,p}$, $CT$-$\mathbb{G}\mathrm{r}_{c,p}$, were already been defined, and we know that there are surcategories too. In particular let us remind that: $T\text{-}\mathbb{G}\mathrm{r}_{c,p} := Pt(T\text{-}\mathbb{G}\mathrm{r}_c)$ and $CT\text{-}\mathbb{G}\mathrm{r}_{c,p} := Pt(CT\text{-}\mathbb{G}\mathrm{r}_c)$.

The purpose of this paragraph is to build the contractibility surmonad. For it we have to demonstrate the following:

**Theorem 11** *The forgetful surfunctor $V$:*

$$(CT\text{-}\mathbb{G}r_{c,p}, \mathtt{A}) \xrightarrow{V} (T\text{-}\mathbb{G}r_{c,p}, \mathtt{A})$$

*has a left suradjoint $H$: $H \dashv V$. Moreover the surfunctor $V$ is surmonadic.* □

To prove this theorem we first have to show its unpointed version:

**Theorem 12** *The forgetful surfunctor $V'$*

$$(CT\text{-}\mathbb{G}r_c, \mathtt{A}) \xrightarrow{V'} (T\text{-}\mathbb{G}r_c, \mathtt{A})$$

*has a left suradjoint $H'$: $H' \dashv V'$.* □

PROOF Let $(P,d,c) \in T\text{-}\mathbb{G}\mathrm{r}_c(0)$. We are going to build the free contractible object $H'((P,d,c)) = (\widetilde{P},\widetilde{d},\widetilde{c};([,]_n)_{n\in\mathbb{N}^*})$ by using elementary techniques resulting from logic. The language that we are going to use does not include variables and will be essentially given by constants (the cells of $P$) and symbols of binary operations noted $[,]_n$ (for each $n \in \mathbb{N}^*$). This construction will give the arrow $(P,d,c) \xrightarrow{c_P} (\widetilde{P},\widetilde{d},\widetilde{c})$ and we shall show that this arrow is universal. The construction of $(\widetilde{P},\widetilde{d},\widetilde{c})$ is made by induction:

- $\widetilde{P}(0)=P(0)$, $\widetilde{d_0}=d_0$, $\widetilde{c}_0=c_0$, also $\widetilde{P}(1)=P(1)$, $\widetilde{d_1}=d_1$, $\widetilde{c}_1=c_1$ and $\widetilde{P}(1) \underset{\widetilde{t}=t}{\overset{\widetilde{s}=s}{\rightrightarrows}} \widetilde{P}(0)$ .

- Let us suppose that our $T$-graph $\widetilde{P}$ is built up to $n$-cells, we are going to show how to build

$$T(G)(n+1) \xleftarrow{\widetilde{d}_{n+1}} \widetilde{P}(n+1) \xrightarrow{\widetilde{c}_{n+1}} G(n+1)$$

At first we can suppose that we have the following $T$-graph of dimension $n$

$$\begin{array}{ccccc}
T(G)(n) & \xleftarrow{\widetilde{d}_n} & \widetilde{P}(n) & \xrightarrow{\widetilde{c}_n} & G(n) \\
s\downarrow\downarrow t & & \widetilde{s}\downarrow\downarrow\widetilde{t} & & s\downarrow\downarrow t \\
T(G)(n-1) & \xleftarrow{\widetilde{d}_{n-1}} & \widetilde{P}(n-1) & \xrightarrow{\widetilde{c}_{n-1}} & G(n-1) \\
\vdots & & \vdots & & \vdots \\
T(G)(1) & \xleftarrow{\widetilde{d}_1} & \widetilde{P}(1) & \xrightarrow{\widetilde{c}_1} & G(1) \\
s\downarrow\downarrow t & & \widetilde{s}\downarrow\downarrow\widetilde{t} & & s\downarrow\downarrow t \\
T(G)(0) & \xleftarrow{\widetilde{d}_0} & \widetilde{P}(0) & \xrightarrow{\widetilde{c}_O} & G(0)
\end{array}$$

In this case we put

$$P(n+1) \subset \widetilde{P}(n+1),\ \widetilde{d}_{n+1}\mid_{P(n+1)} = d_{n+1}, \widetilde{c}_{n+1}\mid_{P(n+1)} = c_{n+1}.$$

Lets note:

$$D_n = \{(\alpha,\beta) \in \widetilde{P}(n) \times \widetilde{P}(n) / (\widetilde{s} \times \widetilde{t})(\alpha,\alpha) = (\widetilde{s} \times \widetilde{t})(\beta,\beta), \text{ and } \widetilde{d_n}(\alpha) = \widetilde{d_n}(\beta)\}.$$

If $(\alpha,\beta) \in D_n$ then $[\alpha,\beta]_n$ is an $(n+1)$-cell (formal) of $\widetilde{P}$ such as [1]:

- $\widetilde{s}([\alpha,\beta]_n) = \alpha$, $\widetilde{t}([\alpha,\beta]_n) = \beta$,
- $\widetilde{d}_{n+1}([\alpha,\beta]_n) = 1_{\widetilde{d_n}(\alpha)=\widetilde{d_n}(\beta)}$, and $\widetilde{c}_{n+1}([\alpha,\beta]_n) = \widetilde{c}_n(\alpha) = \widetilde{c}_n(\beta)$.

Thus by definition $\widetilde{P}(n+1)$ is built with $P(n+1)$ and with all these formal $(n+1)$-cells $[\alpha,\beta]_n$. Also $\widetilde{P}(n+1) \underset{\widetilde{t}}{\overset{\widetilde{s}}{\rightrightarrows}} \widetilde{P}(n)$ will be such as $\widetilde{s}\,|_{P(n+1)} = s$, and $\widetilde{t}\,|_{P(n+1)} = t$.

- We easily show that we have the equalities $\widetilde{d_n}\widetilde{s} = s\widetilde{d}_{n+1}$, $\widetilde{d_n}\widetilde{t} = t\widetilde{d}_{n+1}$, $\widetilde{c}_n\widetilde{s} = s\widetilde{c}_{n+1}$, $\widetilde{c}_n\widetilde{t} = t\widetilde{c}_{n+1}$.

Therefore we obtain a contractible $T$-graph:

$$H'((P,d,c)) := (\widetilde{P}, \widetilde{d}, \widetilde{c}; ([,]_n)_{n\in\mathbb{N}^*}).$$

We also have the canonical embedding $c_P$: $(P,d,c) \overset{(c_P,1_G)}{\hookrightarrow} (\widetilde{P}, \widetilde{d}, \widetilde{c})$ , (i.e $\forall n \in \mathbb{N}, \forall x \in P(n), c_P(x) = x$). Let us prove that $c_P$ is universal.

Let $(Q,d',c';([,]'_n)_{n\in\mathbb{N}^*})$ be another contractible $T$-graph. Let $(P,d,c) \xrightarrow{(f,h)} (Q,d',c')$ be a morphism of $T$-graphs. It means in particular that we give ourselves equalities (in $\infty$-$\mathbb{G}$r) $T(h)d = d'f$, $hc = c'f$. We are going to show that there is a unique morphism of contractible $T$-graphs $(g,h)$ such as $(g,h)(c_P,1_G) = (f,h)$.

**Remark 8** The morphism $h$ is imposed by the way morphisms of $T$-graphs are composed of. Thus the only difficulty is to find $g$. □

[1] Do not forget that we work with an arity $\infty$-graph which is constant.

We are going to build $g$ by induction on $g_n$ which consists of it:

- Because of the definition of $\widetilde{P}(0)$ and of $\widetilde{P}(1)$ we inevitably have $g_0 = f_0$ and $g_1 = f_1$.

- Let us suppose that every $g_i$ are built for $i \in [\![0, n]\!]$. We are going to show that $g_{n+1}$ is then imposed on us. We have the equality $g_{n+1}c_P(n+1) = f_{n+1}$ which immediately shows that $g_{n+1} \mid_{P(n+1)} = f_{n+1}$. As $g$ also preserves contractions we have: If $x = [\alpha, \beta]_n \in \widetilde{P}(n+1) \setminus P(n+1)$, $g_{n+1}(x) := [g_n(\alpha), g_n(\beta)]_n$. Therefore $\forall n \in \mathbb{N}$, $g_n$ does exist and it is unique. Thus we have the existence and the unicity of $g$.

A simple induction allows to show that we have the equalities $T(h)\widetilde{d} = d'g$ and $h\widetilde{c} = c'g$ (in $\infty$-$\mathbb{G}$r). Finaly $c_P$ is universal, what concludes the proof of the theorem 12. ■

Now we are able of showing the main theorem of this paragraph,i.e the pointed version of the theorem which we have just proved. In this demonstration we shall not respect, stricto sansu, the notations of morphisms because there is no risk of confusion.

PROOF Let $(P, d, c; x) \in T\text{-}\mathbb{G}\mathrm{r}_{c,p}(0)$. Let us write

$$H(P, d, c; x) := (H'(P, d, c), c_P \circ x)$$

and let us show that the following arrow (given by the previous arrow $c_P$)

$$(P, d, c; x) \xrightarrow{c_P} H(P, d, c; x)$$

is universal. Let be the arrow $(P, d, c; x) \xrightarrow{f} (Q, d', c'; z, ([,]'_n)_{n \in \mathbb{N}^*})$ such as $(Q, d', c'; z, ([,]'_n)_{n \in \mathbb{N}^*})$ is a contractible pointed $T$-graph. So we are going to prove that there is a unique morphism $H(P, d, c; x) \xrightarrow{g} (Q, d', c'; z, ([,]'_n)_{(n \in \mathbb{N}^*)})$ of $CT\text{-}\mathbb{G}\mathrm{r}_{c,p}$ such as $gc_P = f$. First of all, we are going to show the uniqueness of $g$, which will almost show its existence. If a morphism of $CT\text{-}\mathbb{G}\mathrm{r}_{c,p}$

$$H(P, d, c; x) \xrightarrow{g} (Q, d', c'; z, ([,]'_n)_{n \in \mathbb{N}^*})$$

is such as $gc_P = f$ then in particular the morphism

$$H'(P,d,c) \xrightarrow{g} (Q,d',c';([,]'_n)_{n\in\mathbb{N}^*})$$

of the underlying $T$-graphs, must preserve contractions. We also have the equality $gc_P = f$. But such $g$ exists and is unique, thanks to the universal property of $c_P$. Thus unicity is proved. Let us show the existence of the morphism $H(P,d,c;x) \xrightarrow{g} (Q,d',c';z,([,]'_n)_{n\in\mathbb{N}^*})$ such as $gc_P = f$. Let $H'(P,d,c) \xrightarrow{g} (Q,d',c';([,]'_n)_{n\in\mathbb{N}^*})$ be the previous morphism given by the universal property of $c_P$. In this case we have $I\mathtt{A}(g)1_{I(\mathtt{A}(P))} = I\mathtt{A}(f)$ because $gc_P = f$ and we have $I\mathtt{A}(gc_P) = I\mathtt{A}(g)I\mathtt{A}(c_P) = I\mathtt{A}(g) = I\mathtt{A}(f)$. We also have $gc_Px = zI(\mathtt{A}(g))$ because $gc_Px = fx = zI\mathtt{A}(f) = zI\mathtt{A}(g)$. Thus $H'(P,d,c) \xrightarrow{g} (Q,d',c';([,]'_n)_{n\in\mathbb{N}^*})$ of $\infty$-$\mathbb{G}$r produces the morphism as well

$$H(P,d,c;x) \xrightarrow{g} (Q,d',c';z,([,]'_n)_{n\in\mathbb{N}^*})$$

of $CT$-$\mathbb{G}\mathrm{r}_{c,p}$ such as $gc_P = f$. The surmonadicity of $V$ is an easy application of theorem 4. ∎

The contractibility surmonad will be noted $(C,m,c)$.

Kamel KACHOUR
69 rue Balard, 75015 Paris, France
Phone: (33)(0)1 40 60 95 64
Email:kachour_kamel@yahoo.fr